\RequirePackage{fix-cm}
\documentclass[12pt]{svjour3}                 
\smartqed  
\usepackage{etex} 
\usepackage{graphicx,subfigure,color,colordvi}
\usepackage{latexsym,amsbsy,epsfig,fancybox}
\usepackage{amstext}
\usepackage{amsfonts,textgreek}
\usepackage{mathrsfs}
\usepackage{mathtools}
\usepackage{stmaryrd,dsfont}
\usepackage{bbm}
\usepackage{url}
\usepackage{wrapfig,bm}
\usepackage{tikz}
\usepackage{pstricks}
\usepackage{caption}
\usepackage{psfrag}
\usepackage{upgreek}
\usepackage{isomath}
\usepackage{latexsym}
\usepackage{array}
\usepackage{multirow}
\usepackage{booktabs}
\usepackage{colortbl}
\usepackage{verbatim}
\usepackage[colorlinks=true,citecolor=blue]{hyperref}
\usepackage{epstopdf,overpic}
\usepackage{enumerate}
\usepackage{amsmath,amssymb,xspace}

\usepackage{pgfplots}
\usetikzlibrary{plotmarks}

\usepackage[margin=3cm]{geometry}

\newcommand{\eref}[1]{Equation~(\ref{#1})}
\newcommand{\erefs}[1]{Equations~(\ref{#1})}
\newcommand{\fref}[1]{Figure~\ref{#1}}

\newcommand{\vm}[1]{\bm{\mathrm{#1}}} 
\newcommand{\bsym}[1]{\bm{#1}}

\renewcommand{\Re}{{\rm{I\!R}}}

\newcommand{\mat}[1]{\bm{\mathrm{#1}}} 
\newcommand{\transpose}{\mathrm{T}}

\tikzstyle{nicebox}=[draw=black!100, fill=white!10, rectangle, inner sep=4pt, inner ysep=16pt]
\tikzstyle{niceboxtitle}=[draw=black!100, fill=white, text=black, rectangle]

\usetikzlibrary{arrows,decorations.markings}

\definecolor{forestgreen}{RGB}{34, 139, 34}
\definecolor{lightgray}{gray}{0.92}

\newcommand{\xx}{\vm{x}}
\newcommand{\rmd}{\mathrm{d}}

\newcommand{\cn}{\vm{n}}

\newcommand{\bq}{\mat{q}}
\newcommand{\uu}{\mat{u}}

\begin{document}

\title{A new one point quadrature rule over arbitrary star convex polygon/polyhedron
}

\titlerunning{{One point quadrature over polytopes}}        



\author{S.~Natarajan$^a$, A.~Francis$^a$, E.~Atroshchenko$^b$ S.~P.~A.~Bordas$^{c,d,e}$}

\authorrunning{S.~Natarajan \textit{et al.,}} 

\institute{S.~Natarajan$^a$ \at
              $^a$Department of Mechanical Engineering, Indian Institute of Technology, Madras, Chennai - 600036.\\
              \email{snatarajan@cardiffalumni.org.uk; snatarajan@iitm.ac.in.}           
            \and A.~Francis$^a$ \at
            $^a$ Department of Mechanical Engineering, Indian Institute of Technology, Madras, Chennai - 600036.\\
           \and E.~Atroshchenko$^b$ \at
           $^b$ Department of Mechanical Engineering, University of Chile, Santiago, Chile.\\
           \and
           S.~P.~A.~Bordas$^{c,d,e}$\at
              $^c$Institute of Research and Development, Duy Tan University, K7/25 Quang Trung, Danang, Vietnam.\\
              $^d$Institute of Computational Engineering, Faculty of Science Technology and Communication, University of Luxembourg, Luxembourg. \\ 
}

\date{Received: date / Accepted: date}

\maketitle

\begin{abstract}
The Linear Smoothing (LS) scheme \cite{francisa.ortiz-bernardin2017} ameliorates linear and quadratic approximations over convex polytopes by employing a three-point integration scheme. In this work, we propose a linearly consistent one point integration scheme which possesses the properties of the LS scheme with three integration points but requires one third of the integration computational time. The  essence of the proposed technique is to approximate the strain by the smoothed nodal derivatives that are determined by the discrete form of the divergence theorem. This is done by the Taylor's expansion of the weak form which facilitates the evaluation of the smoothed nodal derivatives acting as  stabilization terms. The smoothed nodal derivatives are evaluated only at the centroid of each integration cell. These integration cells are the simplex subcells (triangle/tetrahedron in two and three dimensions) obtained by subdividing the polytope. The salient feature of the proposed technique is that it requires only $n$ integrations for an $n-$ sided polytope as opposed to $3n$ in~\cite{francisa.ortiz-bernardin2017} and $13n$ integration points in the conventional approach. The convergence properties, the accuracy, and the efficacy of the LS with one point integration scheme are discussed by solving few benchmark problems in elastostatics. 

\keywords{Polygonal finite element method, Wachspress shape functions, numerical integration, linear consistency, one point integration.}
\end{abstract}

\section{Introduction}
\label{intro}
Some of the constraints imposed by the conventional finite element method (FEM) is relaxed by the introduction of elements with arbitrary edges/faces. Approximations on arbitrary polytopes have fueled  the development of Polygonal/Polyhedral Finite Element Methods (PFEM) \cite{sukumar2004,sukumarmalsch2006,talischipaulino2012,Rand_Gillette_Bajaj2014}. POLY elements offer added flexibility in meshing complex geometries through various meshing algorithms using Voronoi tessellation \cite{botschmario2007,talischipereira2014}. Such approaches were used to model complex geometries with inclusions \cite{jaskowiecplucinski2016}, modeling of polycrystalline materials \cite{szesheng2005,botschmario2007,jayabalmenzel2011}. 

Adaptive mesh generation and regeneration such as local refinement and coarsening is also simplified with polytopes, since they naturally address the issues associated with hanging nodes \cite{natarajanooi2015,hoangxu2016}. This has led researchers to develop methods with polygonal discretizations, for example, mimetic finite differences~\cite{lipnikovmanzini2014}, virtual element method~\cite{veigamanzini2012,veigabrezzi2013,veigabrezzi2014}, finite volume method~\cite{droniou2010}, discontinuous Galerkin method~\cite{cangianigeorgoulis2014}, virtual node method~\cite{tangwu2009} and the scaled boundary finite element method~\cite{natarajanooi2014,ooisong2016,natarajanooi2017}. Furthermore, polygonal/polyhedral elements have also been used to solve problems involving large deformations \cite{biabanakikhoei2012}, incompressibility \cite{talischipereira2014}, contact problems \cite{biabanakikhoei2014} and fracture mechanics \cite{khoeiyasbolaghi2015}.

The flexibility provided by polytopes comes with  challenges. First, the arbitrary polytopes usually rely on rational basis functions, i.e., the ratio of two polynomials. The construction of approximation functions over arbitrary polytopes is not unique. These approaches include: Mean value coordinates \cite{floater2003}, Harmonic shape functions \cite{bishop2013}, Laplace basis functions \cite{sukumartabarrei2004} and maximum entropy basis functions \cite{sukumar2004}. 

Integrating such rational functions exactly is not possible in general. One approach to integrate over arbitrary polytopes is to sub-divide the region into triangles (in two dimensions) or tetrahedra (in three dimensions) and then employ conventional quadrature schemes. Although, the approach is simple, it requires \emph{``many''}  integration points to integrate even simple functions~\cite{sukumartabarrei2004,sukumartabarrei2004a,talischipaulino2014}. Moreover, the associated approximations do not pass the patch test ~\cite{sukumartabarrei2004a,talischipaulino2014}.

Inspired by the smoothing technique originally proposed for meshfree methods~\cite{chenwu2001}, a smoothing technique was proposed for polygonal elements in~\cite{dailiu2007,natarajanooi2014}. However, it was shown that the direct application of the smoothing technique with average shape functions does not pass the patch test either and yields less accurate results~\cite{natarajanooi2014}.

Within the framework of the smoothed finite element method (SFEM), Francis \emph{et al.,}~\cite{francisa.ortiz-bernardin2017} proposed a linear smoothing (LS) technique. The LS scheme employs $3n$ and $4n$ integration points for two and three dimensional elements, where $n$ is the number of vertices of the polytope. It was shown with the help of numerical examples that the LS scheme leads to improved accuracy and recovers optimal convergence for the arbitrary convex polytopes. Moreover, it also passes the patch test to machine precision. 

In this paper, we present a new one point quadrature rule over arbitrary star convex polytopes which can reproduce linear strain. In order to achieve this characteristic, the Taylor's expansion of the stiffness matrix and the strain-displacement matrix is employed around the center of the subcell. The modified derivatives are calculated at the centroid of each subcell and a conventional assembly procedure is adopted to calculate the stiffness matrix. The robustness, the accuracy and the convergence properties are studied with a few benchmark problems in elastostatics. The paper is organized as follows: Section \ref{govereqn_weak} presents the governing equations for elasto-statics. Section \ref{oneptdes} presents the new one point quadrature scheme for star convex arbitrary polytopes. Numerical results are presented in Section \ref{numex}, followed by concluding remarks in the last section.

\section{Governing equations for homogeneous linear elastic material}
\label{govereqn_weak}

\subsection{Strong form}

Consider a homogeneous isotropic linear elastic body occupying $d=$ {$2,3$} dimensional space defined by an open domain $\Omega \subset \Re^d$, bounded by the ($d\!-\!1$) dimensional surface $\Gamma$ such that $\Gamma=\Gamma_u \cup \Gamma_t$ and $\emptyset=\Gamma_u\cap\Gamma_t$, where $\Gamma_u$ and $\Gamma_t$ are part of the boundary where Dirichlet and Neumann boundary conditions are specified, with $\vm{n}$ the unit outward normal. The boundary-value problem for linear elastostatics is defined by 
\begin{equation}
\quad \bsym{\nabla} \cdot \bsym{\sigma} + \vm{b} = \textbf{0} ~ \rm{in}~ \Omega,
\label{eqn:problem_strong_form}
\end{equation}
with the following boundary conditions
\begin{align}
\quad \vm{u} &= \bar{\vm{u}} ~ \rm{on}~ \Gamma_u, \nonumber \\
\quad \bsym{\sigma} \cdot \vm{n}
&= \bar{\vm{t}} ~ \rm{on}~ \Gamma_t,
\end{align}
where $\bsym{\sigma}$ is the Cauchy stress tensor and $\vm{u} : \Omega \rightarrow \Re^d$ is the nodal displacement field of the elastic body when it is subjected to external tractions
$\bar{\vm{t}}:\Gamma_t\rightarrow \Re^d$ and body forces $\vm{b}:\Omega\rightarrow\Re^d$. 

\subsection{Weak form}

We first define the infinite dimensional trial ($\mathscr{U}$) and test spaces ($\mathscr{V}$). Let $\mathcal{W}(\Omega)$ be the space including linear displacement fields. 
\begin{subequations}
\begin{align}
\mathscr{U} &:=
\left\{\vm{u}\in [ C^0(\Omega)]^d : \vm{u} \in
[ \mathcal{W}(\Omega)]^d \subseteq [ H^{1}(\Omega)]^d, \ \vm{u} = \bar{\vm{u}}
\ \textrm{on } \Gamma_u \right\},\nonumber \\
\mathscr{V}^0 &:= \left\{\vm{v}\in
[ C^0(\Omega) ]^d : \vm{v} \in [ \mathcal{W}(\Omega)]^d \subseteq
[ H^{1}(\Omega) ]^d, \ \vm{v} = \vm{0} \ \textrm{on } \Gamma_u
\right\}.\nonumber
\end{align}
\end{subequations}
The Bubnov-Galerkin weak form is obtained by testing the strong form \eref{eqn:problem_strong_form} with the test functions in $\mathscr{V}^0$ and integrating over $\Omega$. Using the divergence theorem and the fact that the test functions vanish on the Dirichlet boundary $\Gamma_u$, we obtain the weak form: 
\begin{subequations}
\begin{align}
\begin{split}
\text{Find} \quad \vm{u} \in \mathscr{U} \quad \text{such that, for all} \quad \vm{v} \in \mathscr{V},
\quad a(\vm{u},\vm{v})&= \ell(\vm{v}),
\end{split}\\
\begin{split}
\quad a(\vm{u},\vm{v})&=\int_{\Omega}\bsym{\sigma}(\vm{u}):\bsym{\varepsilon}(\vm{v})\,\rmd V,
\end{split}\\
\begin{split}
\quad \ell(\vm{v})&= \int_{\Omega}\vm{b}\cdot\vm{v}\,\rmd V + \int_{\Gamma_t}\hat{\vm{t}}\cdot\vm{v}\,\rmd S,
\end{split}
\end{align}
\label{eq:weakform}
\end{subequations}
where $\bsym{\varepsilon}=\frac{1}{2} \left[ \nabla \uu + \nabla \uu^{\rm T} \right]$ is the small strain tensor. 

\subsection{Discretisation}

\subsubsection{Discretised weak form}

The domain is partitioned into $n_{el}$ non-overlapping polyhedral elements $\Omega^h$ with planar faces. We define the discrete trial and test spaces by constructing shape functions over the union of all $n_{el}$ $\in$ $\Omega^h$. These shape functions $\phi_e$ are used to discretise the trial and test functions. These trial and test functions are written as a linear combination, over the union of all elements, of the shape functions $\phi_e$ with (vector) coefficients $\vm{u}_e$ :
\begin{align}
\vm{u}^h=\sum_{e=1}^{n_{el}} \phi_e \vm{u}_e \nonumber \\
\vm{v}^h=\sum_{e=1}^{n_{el}} \phi_e \vm{v}_e
\label{eqn:dispfeapprox}
\end{align}
The construction of these (Wachspress) shape functions $\phi_e$ is detailed in Section \ref{shapefunctiondetail}.

With these notations, the following discrete weak form is obtained, which consists in finding $\vm{u}^h \in \mathscr{U}^h \subset \mathscr{U}$ such that for all discretised test functions $\vm{v}^h$ vanishing on the Dirichlet boundary  (in set $\mathscr{V}^{0h} \subset \mathscr{V}^0$),
\begin{equation}\label{eq:discweakform}
 \quad a(\vm{u}^h,\vm{v}^h)= \ell(\vm{v}^h) 
\end{equation}
which leads to the following system of linear equations:
\begin{align}
\mat{K}\mat{u}&=\mat{f},\\
\mat{K}&=\sum_h\mat{K}^h=\sum_h\int_{\Omega^h}\mat{B}^\transpose\mat{C}\mat{B}\,\rmd V, \nonumber\\
\mat{f} &= \sum_h \mat{f}^h=\sum_h\left(\int_{\Omega^h}{\boldsymbol{\phi}}^\transpose\vm{b}\,\rmd V + \int_{\Gamma_t^h}{\boldsymbol{\phi}}^\transpose\hat{\vm{t}}\,\rmd S\right), \nonumber
\label{eq:weakform_disc}
\end{align}
where $\mat{K}$ is the global stiffness matrix, $\mat{f}$ is the global nodal force vector, $\mat{C}$ is the constitutive relation matrix for an isotropic linear elastic material and $\mat{B}=\bsym{\nabla} {\boldsymbol{\phi}}$ is the strain-displacement matrix that is computed using the derivatives of the shape functions. 

\subsubsection{Construction of the shape functions}
\label{shapefunctiondetail}
\begin{figure}[htbp]
\centering
\includegraphics[scale=0.8]{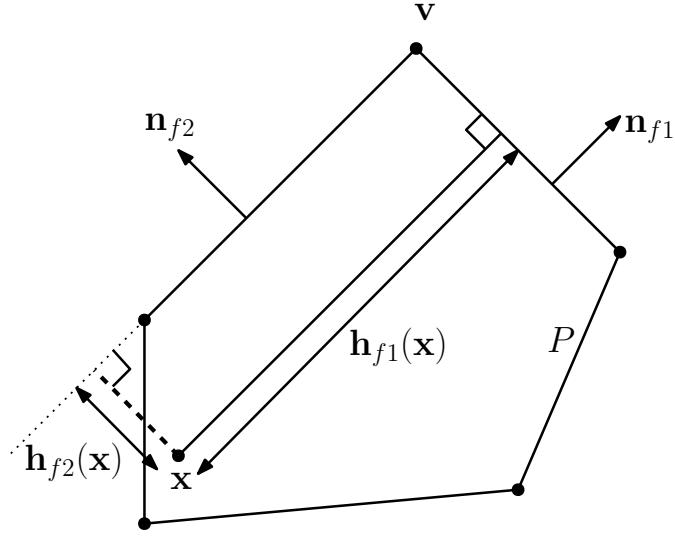}
\caption{Barycentric coordinates: Wachspress basis function}
\label{fig:bary}
\end{figure}
There are different ways to represent the shape functions over arbitrary polytopes~\cite{sukumarmalsch2006}. In this paper, the Wachspress interpolants are used as the approximation functions to describe the unknown fields. These functions are rational polynomials and the construction of the coordinates is as follows: Let $P \subset \Re^3$ be a simple convex polyhedron with facets $F$ and vertices $V$. For each facet $f \in F$, let $\cn_f$ be the unit outward normal and for any $\xx \in P$, let $h_f(\xx)$ denote the perpendicular distance of $\xx$ to $f$, which is given by
\begin{equation}
h_f(\xx) = (\mathbf{v}-\xx) \cdot \cn_f
\end{equation}
for any vertex $\mathbf{v} \in V$ that belongs to $f$. For each vertex $\mathbf{v} \in V$, let $f_1,f_2,f_3$ be the three faces incident to $\mathbf{v}$ and for $\xx \in P$, let
\begin{equation}
w_{\mathbf{v}}(\xx) =  \mathrm{det} ( \mat{p}_{f_1}, \mat{p}_{f_2},\mat{p}_{f_3} )
\end{equation}
where, $\mat{p}_f := \mat{n}_f/h_f(\mat{x})$ is the scaled normal vector, $f_1,f_2,\cdots,f_d$ are the $d$ faces adjacent to $\mat{v}$ listed in an counter-clockwise ordering around $\mat{v}$ as seen from outside $P$ (see \fref{fig:bary}) and $det$ denotes the regular vector determinant in $\mathbb{R}^d$. The shape functions for $\xx \in P$ is then given by
\begin{equation}
\phi_{\mathbf{v}}(\xx) = \frac{ w_{\mathbf{v}}(\xx)}{ \sum\limits_{\mathbf{u} \in V} w_{\mathbf{u}}(\xx)}.
\end{equation}
The Wachspress shape functions are the lowest order shape functions that satisfy boundedness, linearity and linear consistency on convex polytopes~\cite{warren2003,warrenschaefer2007}. On one front, the use of arbitrary shaped elements introduces flexibility and on another, it demands the construction of sufficiently accurate integration rules for computing the terms in the stiffness matrix. This is because the usual and standard integration rules cannot be employed directly. Some of the approaches to integrate over arbitrary polygons include: sub-triangulation~\cite{sukumartabarrei2004a}, Green-Gauss quadrature~\cite{sommarivavianello2007}, nodal quadrature~\cite{Arun2014}, complex mapping~\cite{natarajanbordas2009} conforming interpolant quadrature and strain smoothing~\cite{chenwu2001}. The aforementioned integration rules are restricted to two dimensions. In case of three dimensions, the polyhedron is sub-divided into tetrahedron and cubature rules over the tetrahedron are used for the purpose of numerical integration. Except for the strain smoothing technique, other approaches requires a lot of integration points for sufficient accuracy. In spite of this, it is inferred in~\cite{talischipaulino2014} that the polygonal elements with existing integration technique do not satisfy patch test.

In author's earlier work~\cite{francisa.ortiz-bernardin2017}, a linear smoothing technique was introduced that employed a linear smoothing function and required 3 integration points per subcell in two dimensions and four integration points per subcell in three dimensions. This is accompanied by a modified version of the strain-displacement matrix used to compute the stiffness matrix. The stiffness matrix, as computed within the framework of the SFEM is:
\begin{equation}
\mat{\tilde{K}}=\sum_h \mat{\tilde{K}}^h=\sum_h\int_{\Omega^h}\tilde{\mat{B}}^\transpose\mat{C}\tilde{\mat{B}}\,\rmd V,
\label{eq:modstiffness}
\end{equation}
The next section describes the new one point integration rule to integrate over the arbitrary polytopes.

\section{One point quadrature scheme} 
\label{oneptdes}
\begin{figure}[htpb]
\centering
\begin{subfigure}
\centering
\includegraphics[width=0.85\textwidth]{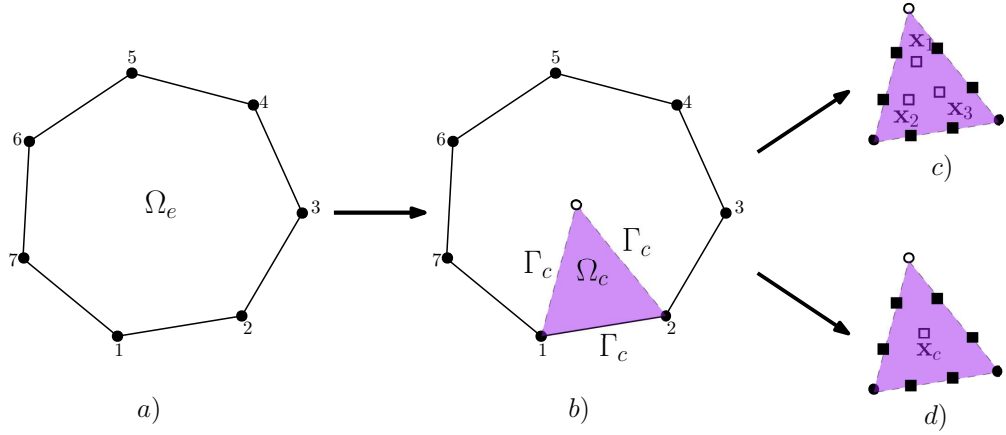}
\caption*{(a) Discertization of arbitrary polygon into triangular subcell using virtual point shown by 'open' circle. }
\end{subfigure}
\quad
\begin{subfigure}
\centering
\includegraphics[width=0.85\textwidth]{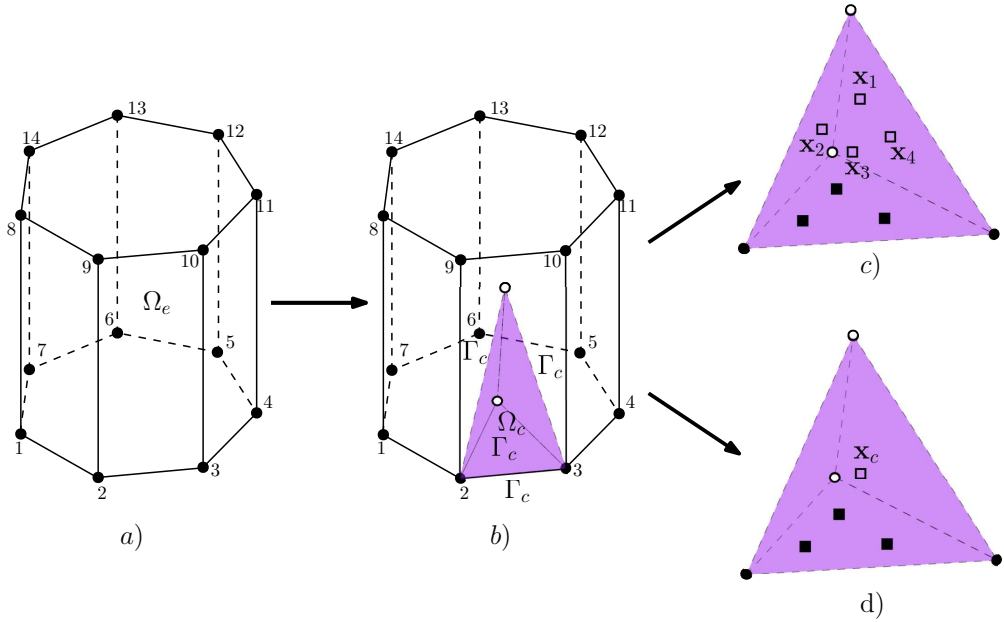}
\caption*{(b) Discertization of arbitrary polyhedron into tetrahedral subcell using virtual point shown by 'open' circle.}
\end{subfigure}\caption{Schematic representation of the three point and one point integration techniques. The nodes are depicted by the filled circles, while the Gauss point per edge/face is shown by filled squares. The smoothed derivatives are computed at the 'open' squares over each smoothing cell denoted by $\Omega_c$.}
\label{fig:subcellRepre}
\end{figure}

In this section, a new numerical integration scheme is proposed to numerically integrate over the star convex arbitrary polygon and polyhedron inspired from the work of Duan \textit{et al.,}~\cite{gaoduan2016}. We restrict ourselves to cell-based smoothing technique, wherein the physical element is sub-divided into simplex elements. This sub-division is solely for the purpose of numerical integration and does not introduce additional degrees of freedom. In this paper, triangles and tetrahedra in two and three dimensions are used as simplex elements. Similar to our earlier work, a linear smoothing function is employed, however, only one integration point is used to compute the modified derivative. This is depicted in \fref{fig:subcellRepre}.
For sake of brevity and simplicity of the notation, the derivation of the proposed smoothing scheme is given in detail only for two-dimensions using the Cartesian coordinate system. The extension to three dimensions is straight forward and interested readers are referred to the corresponding author to obtain the MATLAB code. 

Within the SFEM framework, the discrete modified strain field $\tilde{\varepsilon }_{ij}^h$ that yields the modified strain-displacement matrix $(\tilde{\mat{B}})$ which is then used to build the stiffness matrix is related to the compatible strain field ${\varepsilon }_{ij}^h$ by:
\begin{equation}
\tilde {\varepsilon }_{ij}^h (\xx)=\int_{\Omega_C^h}
{\varepsilon _{ij}^h (\xx) ~\bq(\xx) \rmd V }
\label{eqn:epsilonvar}
\end{equation}
where $q(\xx)$ is the smoothing function. On writing \eref{eqn:epsilonvar} at the basis functions derivative level and invoking Gauss-Ostrogradsky theorem, we get:
\begin{equation}
\int_{\Omega_C^h}\phi_{I,x}~\bq(\vm{x})\, \rmd V  =  \int_{\Gamma_C^h}\phi_I~ \bq(\vm{x})n_j\, \rmd S - \int_{\Omega_C^h}\phi_I ~\bq_{,x}(\vm{x})\, \rmd V
\label{eq:divconsistency}
\end{equation}
In this work, a linear smoothing function $\bq(\xx) = \left\{1,~x,~y \right\}$ in two dimensions and $\bq(\xx) = \left\{1,~x,~y,~z\right\}$ in three dimensions is employed and numerical integration is employed to evaluate the terms in \eref{eq:divconsistency}. Note that the domain integral in \eref{eq:divconsistency} is evaluated at the center of the subcell, $\xx_c = (x_c,y_c)$ (see \fref{fig:subcellRepre}. The center of the subcell is denoted by `open' circle), whilst, the boundary integral is evaluated along the boundary of the subcell (the location of integration point on the boundary is represented by `filled' square in \fref{fig:subcellRepre}). However, this will lead to a singular system~\cite{}. This is circumvented by introducing higher order derivatives, viz., $\tilde{\phi}_{I,x}(x_c), \tilde{\phi}_{I,xx}(x_c), \tilde{\phi}_{I,xy}(x_c)$ by using Taylor's expansion of the modified derivatives around the center of the subcell, $\xx_c$. The Taylor's expansion (around the center of the subcell $\xx_c$) of $\tilde{\phi}_{I,x}(\xx),~\bq(\xx)$ and $\phi_I(\xx)$, used is defined as:
\begin{subequations}
\begin{align}
\tilde{\phi}_{I,x}(\xx)  &= \tilde{\phi}_{I,x}(\xx_c) + (x-x_c) \tilde{\phi}_{I,xx}(\xx_c) + (y-y_c)\tilde{\phi}_{I,xy}(\xx_c) + \mathcal{O}( (\xx-\xx_c)^2) \\ 
\bq(\xx) &= \bq(\xx_c) + (x-x_c) \bq_{,x}(\xx) + (y-y_c)\bq_{,y}(\xx) \\
\begin{split}
\phi_{I}(\xx)  &= \phi_{I}(\xx_c) + (x-x_c) \phi_{I,x}(\xx_c) + (y-y_c)\phi_{I,y}(\xx_c) + \frac{1}{2} (x-x_c)^2 \phi_{I,xx}(\xx_c)  \\  &+  (x-x_c)(y-y_c)\phi_{I,xy}(\xx_c) + \frac{1}{2} (y-y_c)^2 \phi_{I,yy}(\xx_c) + \mathcal{O}( (\xx-\xx_c)^3)
\end{split}
\end{align}
\label{eq:taylorexpansion}
\end{subequations}
Upon substituting \eref{eq:taylorexpansion} into \eref{eq:divconsistency}, we obtain:
\begin{equation}
\begin{split}
\bq(\xx_c)A\tilde{\phi}_{I,x}(\xx_c) + \left[ \bq_{,x}(\xx_c) I_c^{xx} + \bq_{,y}(\xx_c) I_c^{xy} \right] \tilde{\phi}_{I,xx}(\xx_c) + \left[ \bq_{,x}(\xx_c) I_c^{xy} + \bq_{,y}(\xx_c)I_c^{yy} \right] \tilde{\phi}_{I,xy}(\xx_c) \\ = \int\limits_{\Gamma_{c}^h} \phi_I(\xx) \bq(\xx) \mathbf{n}~\mathrm{d}\Gamma - \left[ A\phi_I(\xx_c) + \frac{1}{2}I_c^{xx} \phi_{I,xx}(\xx_c) + I_c^{xy}\phi_{I,xy}(\xx_c) + \frac{1}{2} I_c^{yy} \phi_{I,yy}(\xx_c) \right]
\end{split}
\label{eq:newtaylorform}
\end{equation}
where $A = \int\limits_{\Omega_{c}^h} \mathrm{d}\Omega$ is the area of the integration domain $\Omega$. The first order area moments with respect to cell center $\xx_c$ vanish and the second order area moments are given by: 
\begin{equation}
\left\{ \begin{array}{c} I_c^{xx} \\ I_c^{xy} \\ I_c^{yy} \end{array} \right\} = \int\limits_{\Omega_c} \left\{ \begin{array}{c} (x-x_c)^2 \\ (x-x_c)(y-y_c) \\ (y-y_c)^2 \end{array} \right\}~\mathrm{d}\Omega
\end{equation}
\begin{remark}
For a regular polygon, the second order area moment $I_c^{xy}$ also vanish apart from the first order area moments.
\end{remark}
\begin{remark}
\eref{eq:lineareqn} always have an unique solution provided the triangles do not degenerate to a line. 
\end{remark}
This now leads to the following system of linear equations:
\begin{equation}
\mat{W} \mat{d}_j = \mat{f}_j, \quad j = 1,2
\label{eq:lineareqn}
\end{equation}
where,
\begin{align}
\mat{W} &= \left[ \begin{array}{ccc} A & 0 & 0 \\ Ax_c & I_c^{xx} & I_c^{xy} \\ Ay_c & I_c^{xy} & I_c^{yy} \end{array} \right] \nonumber \\ \nonumber \\
\mat{f}_1 &= \left\{ \begin{array}{c} \sum\limits_{L=1}^3 \sum\limits_{G=1}^2 \phi_I(\xx_G)n_x^Lw_G \\ \sum\limits_{G=1}^2 \phi_I(\xx_G)x_G n_x^Lw_G - Fg \\ \sum\limits_{G=1}^2 \phi_I(\xx_G) y_G n_y^Lw_G \end{array} \right\} \nonumber \\ \nonumber \\
\mat{f}_2 & = \left\{ \begin{array}{c} \sum\limits_{L=1}^3 \sum\limits_{G=1}^2 \phi_I(\xx_G)n_x^Lw_G \\ \sum\limits_{G=1}^2 \phi_I(\xx_G)x_G n_x^Lw_G \\ \sum\limits_{G=1}^2 \phi_I(\xx_G) y_G n_y^Lw_G - Fg \end{array} \right\}
\end{align}
where 
\begin{equation*}
Fg = A\phi_I(\xx_c) + \frac{1}{2}I_c^{xx}\phi_{I,xx}(\xx_c) + I_c^{xy} \phi_{I,xy}(\xx_c) + \frac{1}{2} I_c^{yy} \phi_{I,yy}(\xx_c)
\end{equation*}
where $\phi_I(\xx_c),\, \phi_{I,xx}(\xx_c) \, \phi_{I,yy}(\xx_c)$ and $\phi_{I,xy}(\xx_c)$ are the barycentric coordinates and its derivatives are evaluated at the center of the cell, $(x_G,~ y_G)$  and $w_G$ are the integration points and the weights respectively, along the boundary of the smoothing cells (see ~\fref{fig:subcellRepre}. The integration points are shown as filled circles) and $n_x^L$ and $n_y^L$ are the outward normals along the boundary of the smoothing cell. The solution vector is given by:
\begin{subequations}
\begin{align}
\mat{d}_1 &= \left\{ \begin{array}{c} \tilde{\phi}_{I,x}(\xx_c) \\ \tilde{\phi}_{I,xx}(\xx_c) \\ \tilde{\phi}_{I,xy}(\xx_c) \end{array} \right\} \\
\mat{d}_2 &= \left\{ \begin{array}{c} \tilde{\phi}_{I,y}(\xx_c) \\ \tilde{\phi}_{I,yx}(\xx_c) \\ \tilde{\phi}_{I,yy}(\xx_c) \end{array} \right\}
\end{align}
\end{subequations}
This is further used to construct the modified strain displacement matrix and its derivatives used to evaluate the stiffness matrix as:
\begin{equation}
\tilde{\mat{B}}=\begin{bmatrix}
\tilde{\mat{B}}_1 & \tilde{\mat{B}}_2 &....& \tilde{\mat{B}}_n
\end{bmatrix}
\end{equation}
\begin{equation}
\tilde{\mat{B}}_I(\xx_c)=\begin{bmatrix}
 \tilde{\phi}_{I,x}(\xx_c)& 0\\ 0 & \tilde{\phi}_{I,y}(\xx_c)\\\tilde{\phi}_{I,y}(\xx_c) & \tilde{\phi}_{I,x}(\xx_c)
\end{bmatrix}
\end{equation}
\begin{equation}
\frac{\partial \tilde{\mat{B}}_I(\xx_c)}{\partial x}=\begin{bmatrix}
 \tilde{\phi}_{I,xx}(\xx_c)& 0\\ 0 & \tilde{\phi}_{I,yx}(\xx_c)\\\tilde{\phi}_{I,yx}(\xx_c) & \tilde{\phi}_{I,xx}(\xx_c)
\end{bmatrix}
\end{equation}
\begin{equation}
\frac{\partial \tilde{\mat{B}}_I(\xx_c)}{\partial y}=\begin{bmatrix}
 \tilde{\phi}_{I,xy}(\xx_c)& 0\\ 0 & \tilde{\phi}_{I,yy}(\xx_c)\\\tilde{\phi}_{I,yy}(\xx_c) & \tilde{\phi}_{I,xy}(\xx_c)
\end{bmatrix}
\end{equation}
It should be noted that in the proposed technique the smoothed nodal derivatives are used to compute the terms in the modified stiffness matrix. To introduce the higher order modified derivatives into the final discretized form (see \eref{eq:modstiffness}), the stiffness matrix is expanded in Taylor's series (around the center of the subcell, $\xx_c$) as:
\begin{subequations}\label{eq:taylorsexp}
\begin{align}
\mat{\tilde{K}^{\Omega_c}}&=\int_{\Omega_c}\tilde{\mat{B}}^\transpose\mat{C}\tilde{\mat{B}}\ \rmd V, \label{eq:LS}\\ 
&=\int_{\Omega_c}\left[ \tilde{\mat{B}}^\transpose + \frac{\partial  \tilde{\mat{B}}^\transpose}{\partial x} (x-x_c) + \frac{\partial  \tilde{\mat{B}}^\transpose}{\partial y} (y-y_c) \right] \mat{C} \left[ \tilde{\mat{B}}^\transpose + \frac{\partial  \tilde{\mat{B}}^\transpose}{\partial x} (x-x_c) + \frac{\partial  \tilde{\mat{B}}^\transpose}{\partial y} (y-y_c) \right] \rmd V, \label{eq:onept}
\end{align}
\end{subequations}
However, to compute the body forces, standard Wachspress interpolants and its higher order derivatives are employed. The body force is computed as follows:
\begin{subequations}
\label{eq:taylorsexpforce}
\begin{align}
\mat{f^{b}}&=\int_{\Omega_c} \left( {{\boldsymbol{\phi}}}^\transpose{\mat{b}}\right)  \rmd V, \label{eq:LSb}\\ 
\begin{split}
&=\int_{\Omega_c} \left\{ {\boldsymbol{\phi}}^\transpose{\mat{b}}\vert_{(\xx_c)} + \frac{\partial  {\boldsymbol{\phi}}^\transpose}{\partial x} \mat{b} \vert_{(\xx_c)} (x-x_c)
 + \frac{\partial  {\boldsymbol{\phi}}^\transpose}{\partial y} \mat{b} \vert_{(\xx_c)} (y-y_c) \right.\\
  &+\left. \frac{1}{2} \frac{\partial^2 {\boldsymbol{\phi}}^\transpose}{\partial x^2} \mat{b} \vert_{(\xx_c)} (x-x_c)^2
 +\frac{1}{2} \frac{\partial^2 {\boldsymbol{\phi}}^\transpose}{\partial y^2} \mat{b} \vert_{(\xx_c)} (y-y_c)^2 \right.\\
 &+ \left. \frac{\partial^2 {\boldsymbol{\phi}}^\transpose}{\partial x \partial y} \mat{b} \vert_{(\xx_c)} (x-x_c)(y-y_c) \right\} ~ \rmd V ,
\end{split} \\ 
 &=\mat{A}{\boldsymbol{\phi}}^\transpose {\mat{b}}\vert_{(\xx_c)}+\frac{1}{2}I_c^{xx}\frac{\partial^2 {\boldsymbol{\phi}}^\transpose}{\partial x^2} \mat{b} \vert_{(\xx_c)} + \frac{1}{2} I_c^{yy} \frac{\partial^2 {\boldsymbol{\phi}}^\transpose}{\partial y^2}  \mat{b} \vert_{(\xx_c)} + I_c^{xy} \frac{\partial^2 {\boldsymbol{\phi}}^\transpose}{\partial x \partial y} \mat{b} \vert_{(\xx_c)}
\end{align}
\end{subequations}

\section{Numerical examples}
\label{numex}
In this section, we demonstrate the accuracy and the convergence properties of the proposed linear smoothing scheme (LS) over arbitrary polytopes using 1\textit{n} integration point. The LS scheme is compared to the constant smoothing (CS) scheme by solving few benchmark problems. We also demonstrate the performance of the proposed scheme in a simple three-dimensional elasticity problem. In all the numerical examples, we discretize the domain with arbitrary polytopes based on centroid Voronoi tessellation. The two dimensional polygonal meshes are generated by using the built-in Matlab function {\scriptsize{voronoin}} and the Matlab functions in Polytop~\cite{talischipaulino2012}. The open-source software Neper~\cite{queydawson2011} is used to generate polyhedra meshes. For the purpose of error estimation and convergence studies, the $L^2$ norm and $H^1$ seminorm of the error are used. The following convention is used while discussing the results:
\begin{itemize}
\item CS: constant smoothing over arbitrary polygons in two dimensions.
\item LS3$n$-2D, LS3$n$-3D: linear smoothing scheme with three point integration rule over arbitrary polytopes, in two and in three dimensions, respectively. 
\item LS1-2D, LS1-3D: linear smoothing scheme with one point integration rule over arbitrary polytopes, in two and three dimensions, respectively. 
\end{itemize}

\begin{figure}[htpb]
\centering
\subfigure[]{\includegraphics[width=0.47\textwidth]{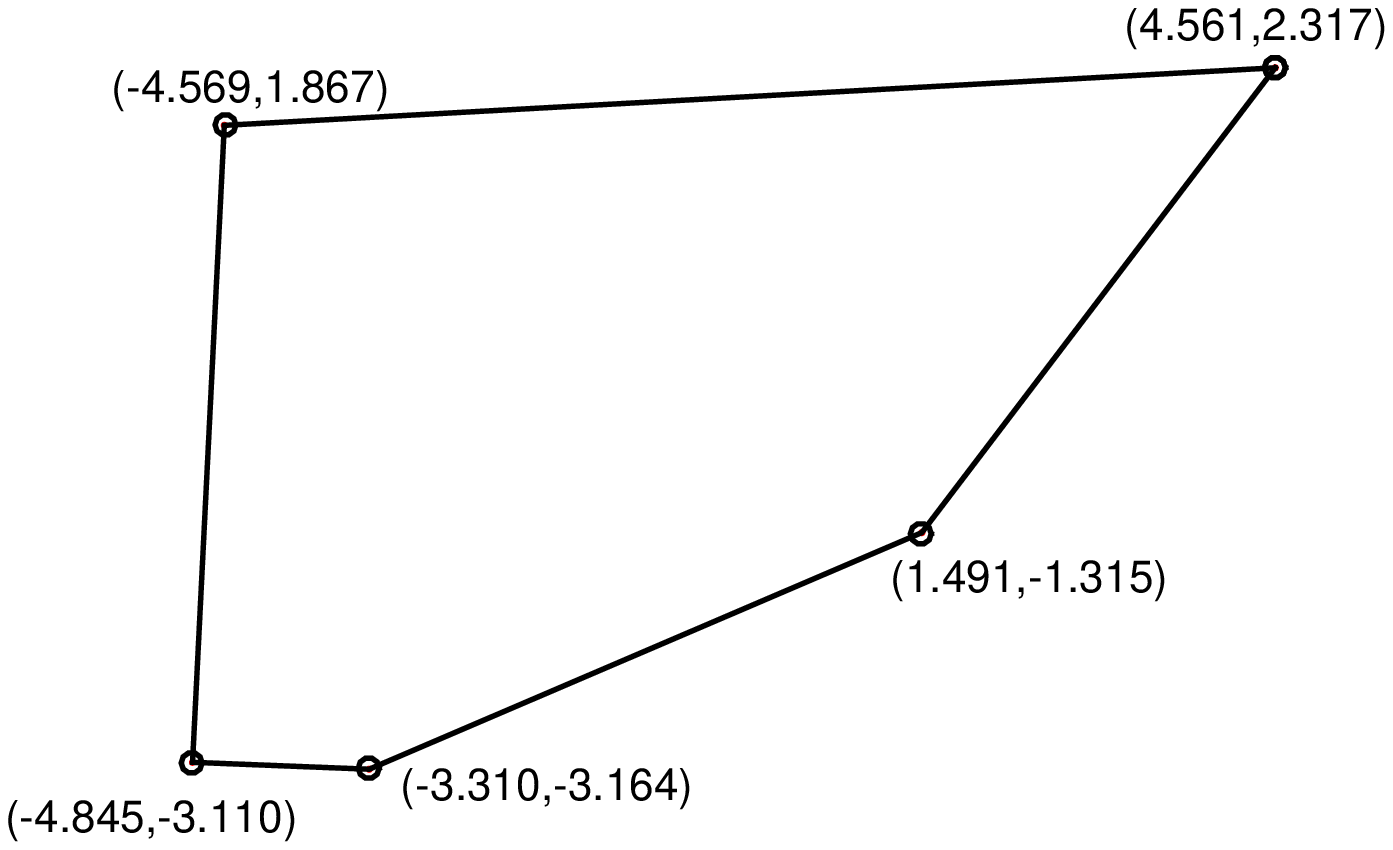}}
\subfigure[]{\includegraphics[width=0.47\textwidth]{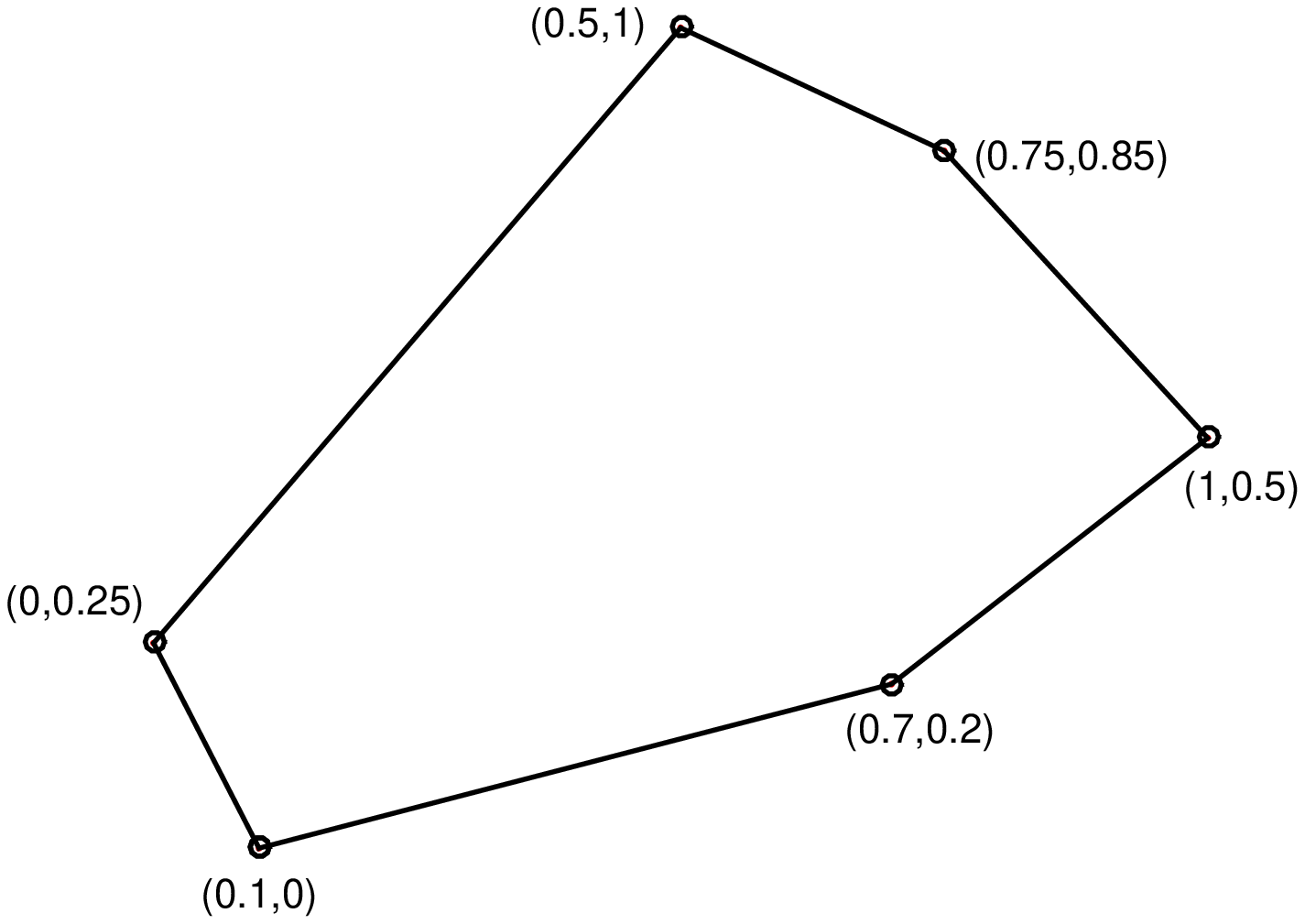}}
\subfigure[]{\includegraphics[width=0.47\textwidth]{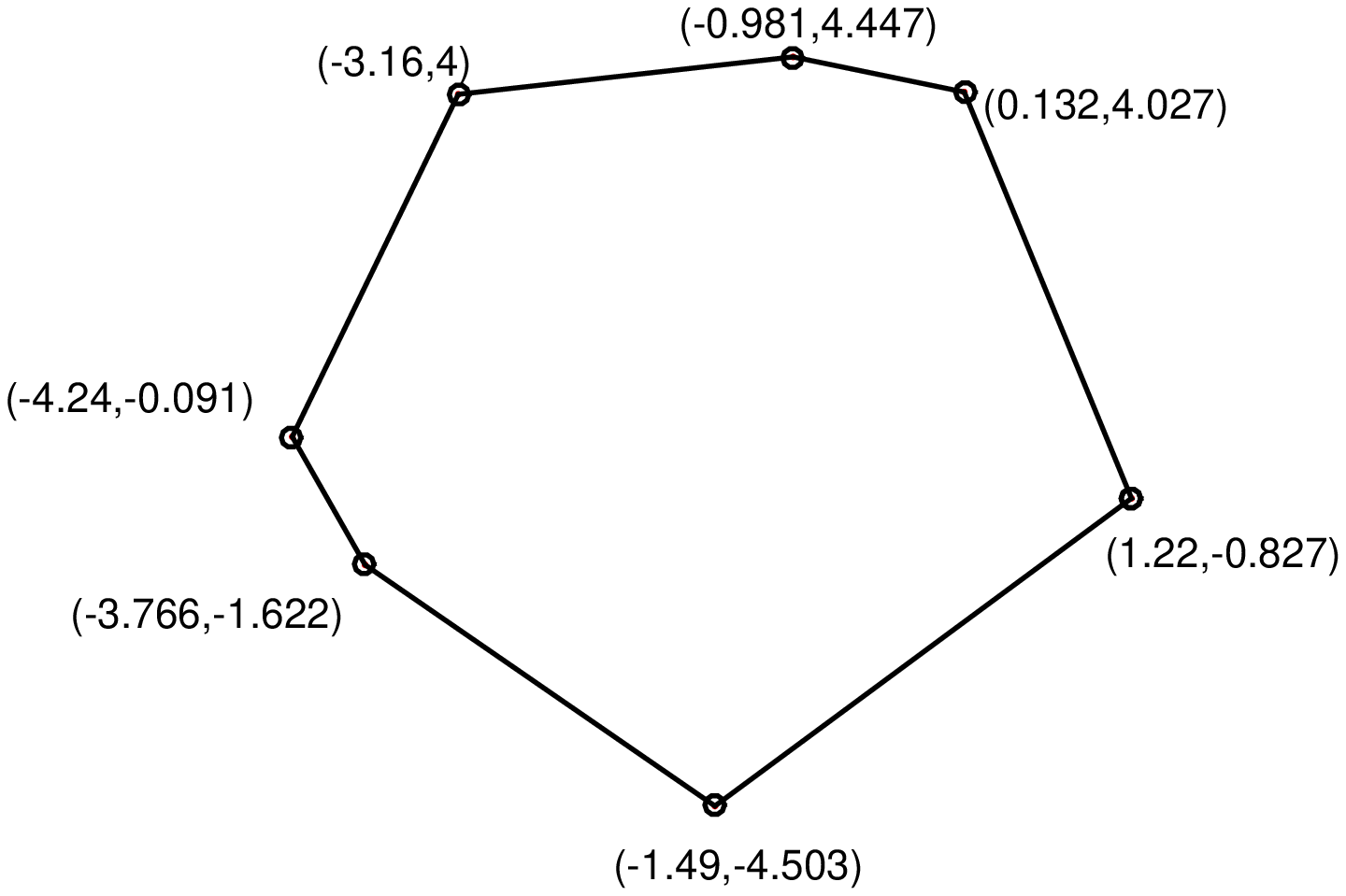}}
\subfigure[]{\includegraphics[width=0.47\textwidth]{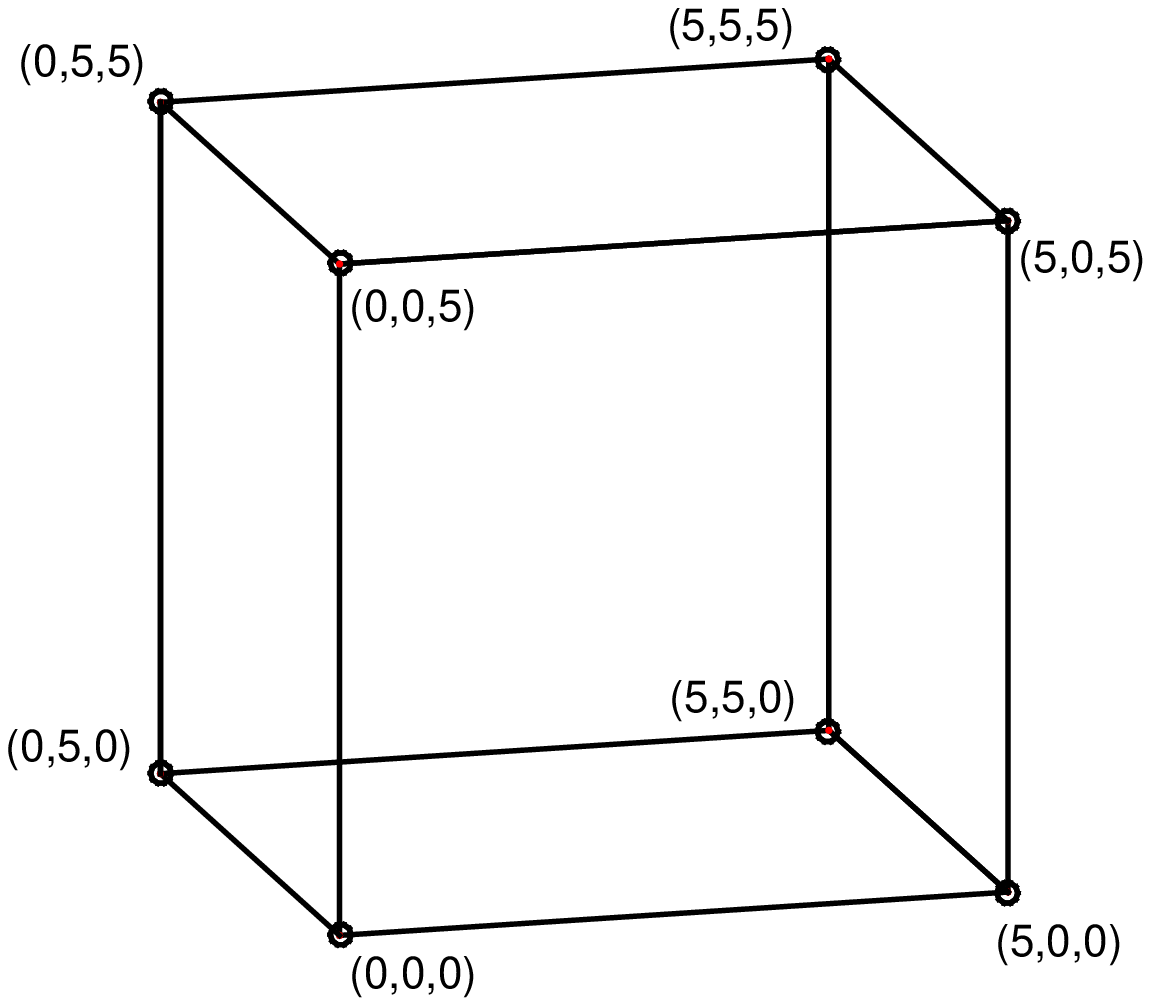}}
\caption{Arbitrary polytopes: a) Pentagon b) Hexagon c) Heptagon and d) Hexahedron.}
\label{fig:arbitrarypoly}
\end{figure}

Before we proceed with the numerical examples, the proposed integration scheme is employed to numerically integrate few polynomials over arbitrary polygons and polyhedra (see \fref{fig:arbitrarypoly} for description of polytopes). The geometry of the polygon and the polynomials are take from~\cite{natarajanbordas2009,chinlasserre2015}. The results from the proposed method are compared with conforming interpolant quadrature (CIQ)~\cite{natarajanooi2017} and analytical solutions. Tables \ref{table:pentagon} - \ref{table:H8} shows the results of numerical integration of the polynomials with the proposed scheme. It is opined that the proposed numerical integration yields accurate results when compared to conventional integration. It is further emphasized that the proposed approach requires only $n$ integration points   

\begin{table}
\renewcommand{\arraystretch}{1.5}
\caption{Numerical integration of polynomial functions over irregular pentagonal domain: comparison between the CIQ using 65 integration points and the LS1-2D using 5 integration points.}
\centering
\begin{tabular}{lcccccc}
\hline
Polynomial & Exact results & \multicolumn{2}{c}{CIQ} &  & \multicolumn{2}{c}{LS1-2D}  \\
\cline{3-4}\cline{6-7}
function & & Evaluated & Relative  && Evaluated & Relative \vspace{-0.15cm}\\
& &results & error &&  results & error   \\
\hline
$1$ & 32.95749050000000 & 32.95749384993172 & 1.02$\times$10$^{-07}$  &&  32.95749050000000 &0.00000 \\
$x$ & 36.57243417200000 & 36.57238829642015 & 1.25$\times$10$^{-06}$ && 36.57243417200428 & 1.17$\times$10$^{-13}$ \\
$x^2$ & 212.9212361315146 & 212.9188256614823& 1.13$\times$10$^{-05}$ && 212.9212361315097&  2.30$\times$10$^{-14}$ \\
$xy$ & 47.43617079993337 & 47.43672095414618 & 1.16$\times$10$^{-05}$  && 47.43617079994323 & 2.08$\times$10$^{-13}$  \\
\hline
\end{tabular}
\label{table:pentagon}
\end{table}

\begin{table}
\renewcommand{\arraystretch}{1.5}
\caption{Hexagonal domain: comparison between the CIQ using 78 integration points and the LS1-2D using 6 integration points.}
\centering
\begin{tabular}{lcccccc}
\hline
Polynomial & Exact results & \multicolumn{2}{c}{CIQ} &  & \multicolumn{2}{c}{LS1-2D}  \\
\cline{3-4}\cline{6-7}
function &~ & Evaluated & Relative  && Evaluated & Relative \vspace{-0.15cm}\\
&~&results & error && results & error  \\
\hline
$1$ & 0.535000000000000 & 0.5349995495917054 & 8.42$\times$10$^{-07}$  && 0.5350000000000000 & 0.00000  \\
$x$ & 0.261416666666667 & 0.2614159504448907& 2.74$\times$10$^{-06}$  && 0.2614166666666675&  1.70$\times$10$^{-15}$  \\
$x^2$ & 0.154606250000000 & 0.1546051808604657 & 6.92$\times$10$^{-06}$  && 0.1546062500000020 & 1.28$\times$10$^{-14}$  \\
$xy$ & 0.133510416666667 & 0.1335091506087918 & 9.48$\times$10$^{-06}$  && 0.1335104166666682 & 9.15$\times$10$^{-15}$  \\
\hline
\end{tabular}
\label{table:hexagon}
\end{table}

\begin{table}
\renewcommand{\arraystretch}{1.5}
\caption{Irregular heptagon domain: comparison between the CIQ 91 integration points and the LS1-2D using 7 integration points.}
\centering
\begin{tabular}{lcccccc}
\hline
Polynomial & Exact results & \multicolumn{2}{c}{CIQ} &  & \multicolumn{2}{c}{LS1-2D}  \\
\cline{3-4}\cline{6-7}
function&~ & Evaluated & Relative  && Evaluated & Relative \vspace{-0.15cm}\\
&~&results & error && results & error  \\
\hline
$1$ & 32.368828500000001 & 32.36886234983820 & 1.05$\times$10$^{-06}$  && 32.36882850000001 & 0.00000  \\
$x$ & 49.533099820500006 & 49.53510843558120 & 4.06$\times$10$^{-05}$  && 49.53309982050202& 4.07$\times$10$^{-14}$  \\
$x^2$ & 126.2695344633893 & 126.2673175334039 & 1.76$\times$10$^{-05}$ && 126.2695344633924 & 2.49$\times$10$^{-14}$  \\
$xy$ & 18.035793954103632 & 18.03503793073484 & 4.19$\times$10$^{-05}$  && 18.03579395410433 & 3.86$\times$10$^{-14}$  \\
\hline
\end{tabular}
\label{table:heptagon}
\end{table}


\begin{table}
\renewcommand{\arraystretch}{1.2}
\caption{Three dimensional hexahedron: comparison between the CIQ and the LS1-3D. The polynomial used for the purpose of integration is: $x^2 + y^2 + xy + z^2$.}
\centering
\begin{tabular}{l p{3.5cm} p{3.5cm} }
 \hline
 ~ & CIQ & LS1-3D\\
 \hline
Exact results & 3906.25 & 3906.25\\
Number of integration points & 324 & 24 \\
Evaluated results & 3906.250000000005 & 3906.250002947967\\
Relative error& 1.16$\times$10$^{-15}$ & 7.55$\times$10$^{-10}$ \\
\hline
\end{tabular}
\label{table:H8}
\end{table}

\subsection{Linear patch test}
In the first example, the accuracy and the convergence properties of the proposed one point quadrature (LS1-2D, LS1-3D) is demonstrated with a linear and a quadratic patch test.
\paragraph{Linear patch test}
The following displacements are prescribed on the boundary in the two-dimensional case:
\begin{equation}
\begin{pmatrix} \hat{u} \\ \hat{v} \end{pmatrix} = \begin{pmatrix} 0.1+0.1x+0.2y \\ 0.05+0.15x + 0.1y \end{pmatrix}
\end{equation}
and in the three-dimensional case the following displacements are prescribed on the boundary:
\begin{equation}
\begin{pmatrix} \hat{u} \\ \hat{v} \\ \hat{w} \end{pmatrix} = \begin{pmatrix} 0.1+0.1x+0.2y+0.2z \\ 0.05+0.15x+0.1y+0.2z \\ 0.05+0.1x+0.2y+0.2z \end{pmatrix}.
\end{equation}

\begin{figure}
\centering
\includegraphics[width=1\textwidth]{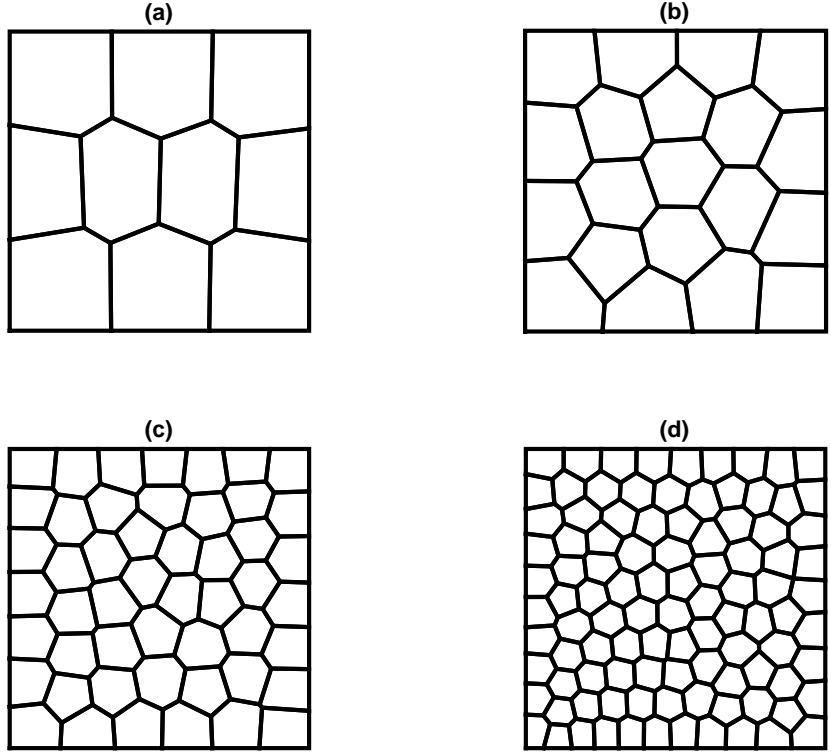}
\caption{Square domain discretized with polygonal elements. Representative meshes containing (a) 10,  (b) 20, (c) 50 and (d) 100 polygons.}
\label{fig:pmesh2d}
\end{figure}

\begin{figure}
\centering
\subfigure[]{\includegraphics[scale=0.17]{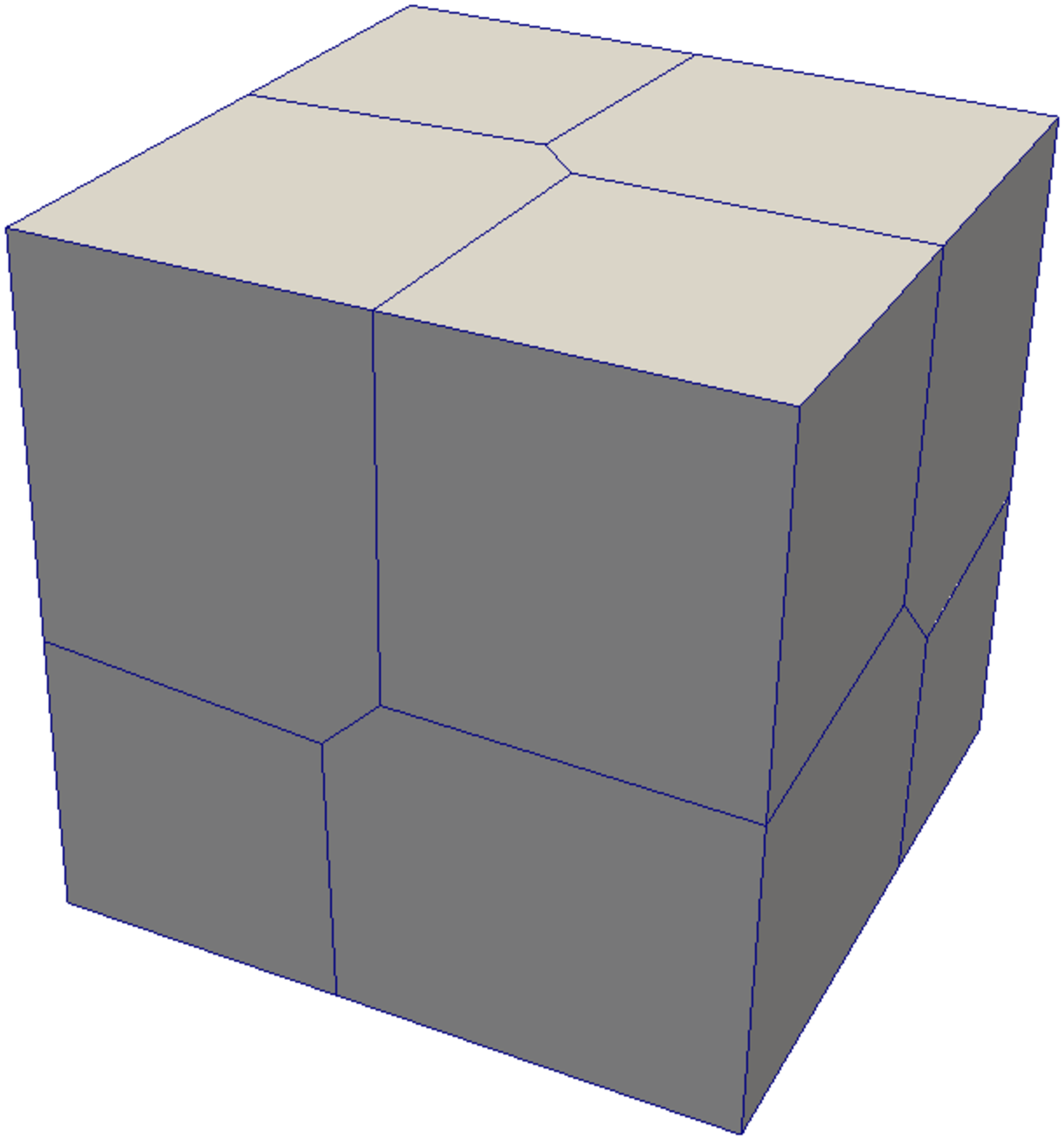}}
\subfigure[]{\includegraphics[scale=0.17]{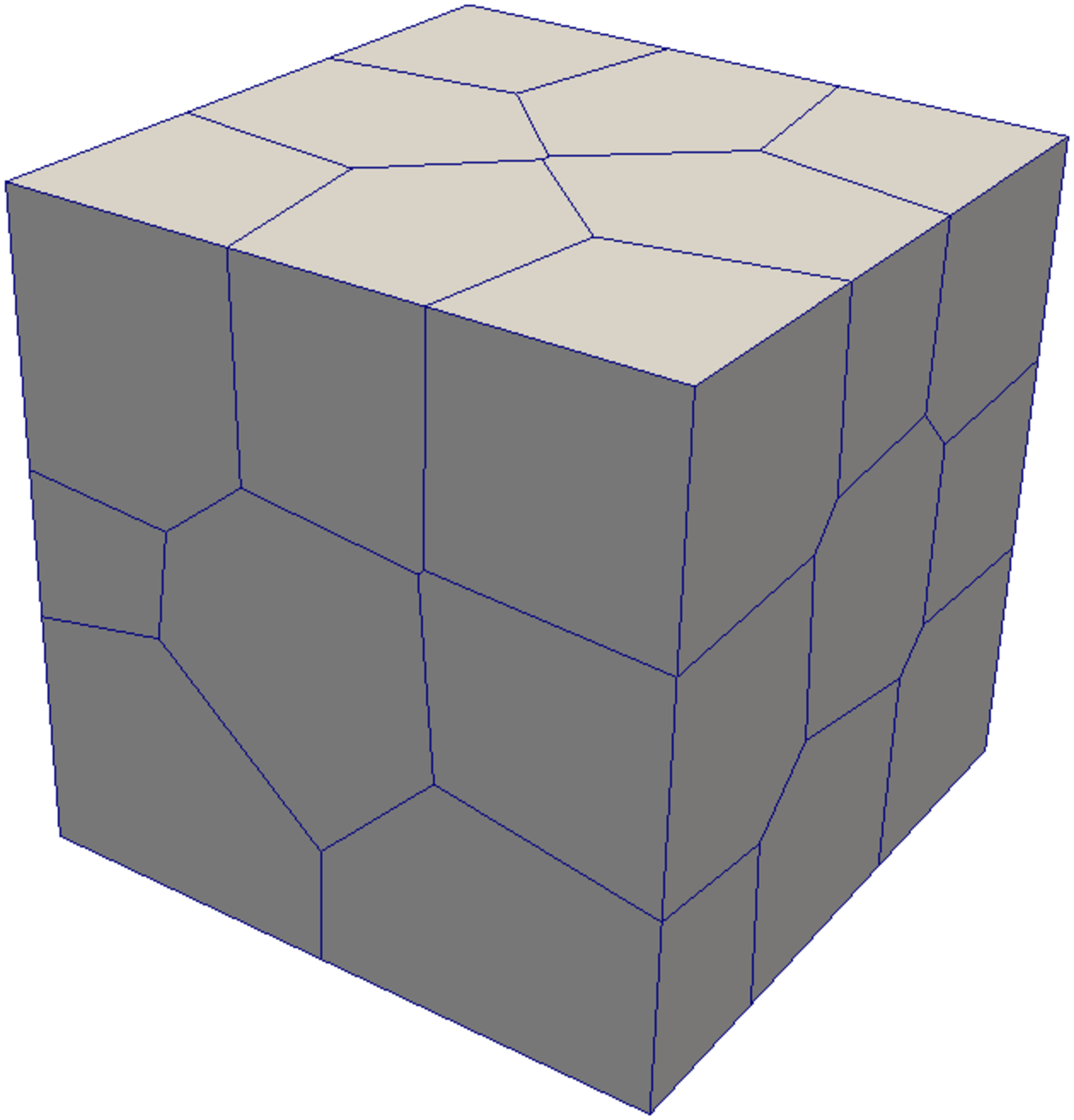}}
\subfigure[]{\includegraphics[scale=0.17]{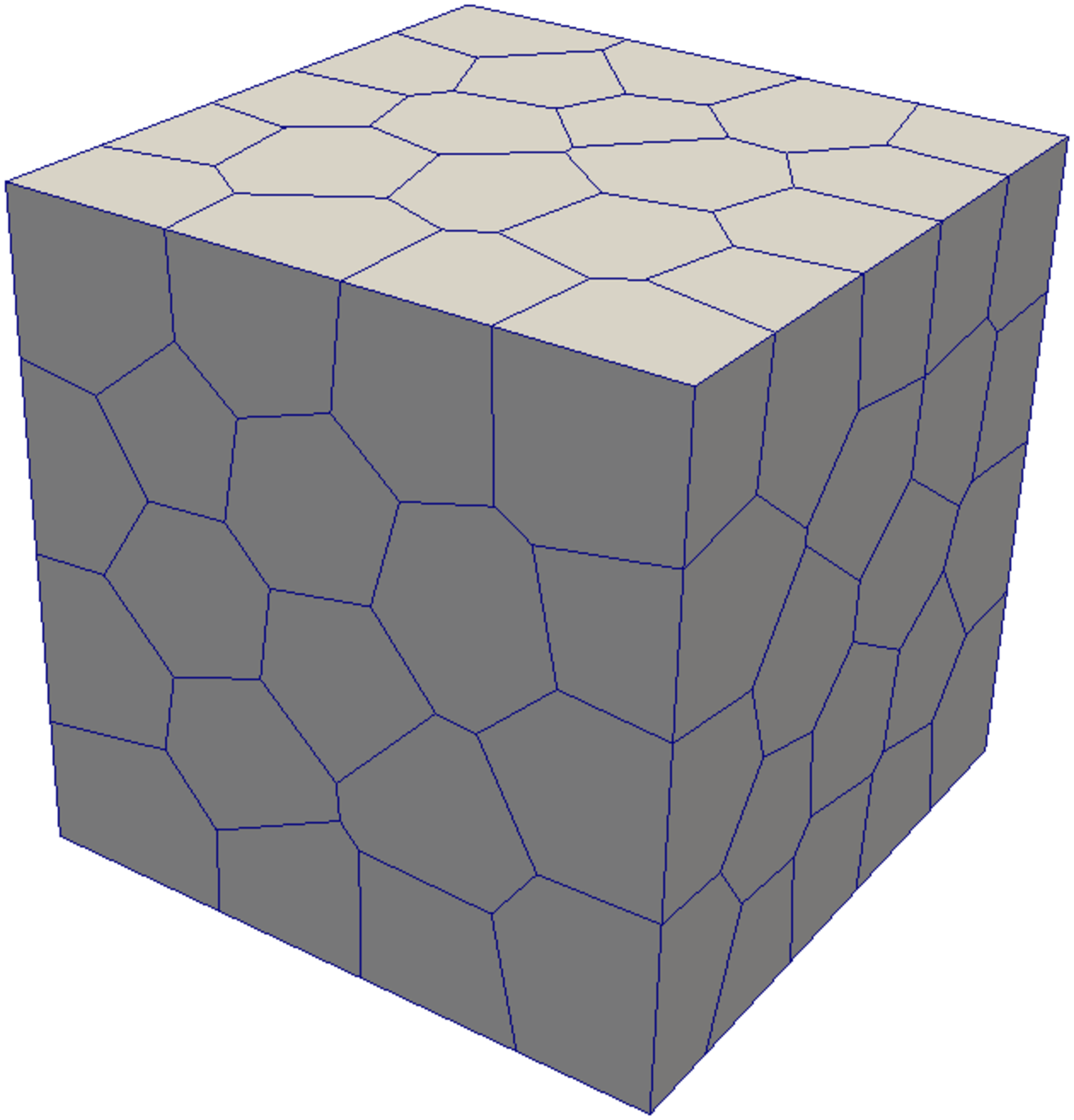}}
\subfigure[]{\includegraphics[scale=0.17]{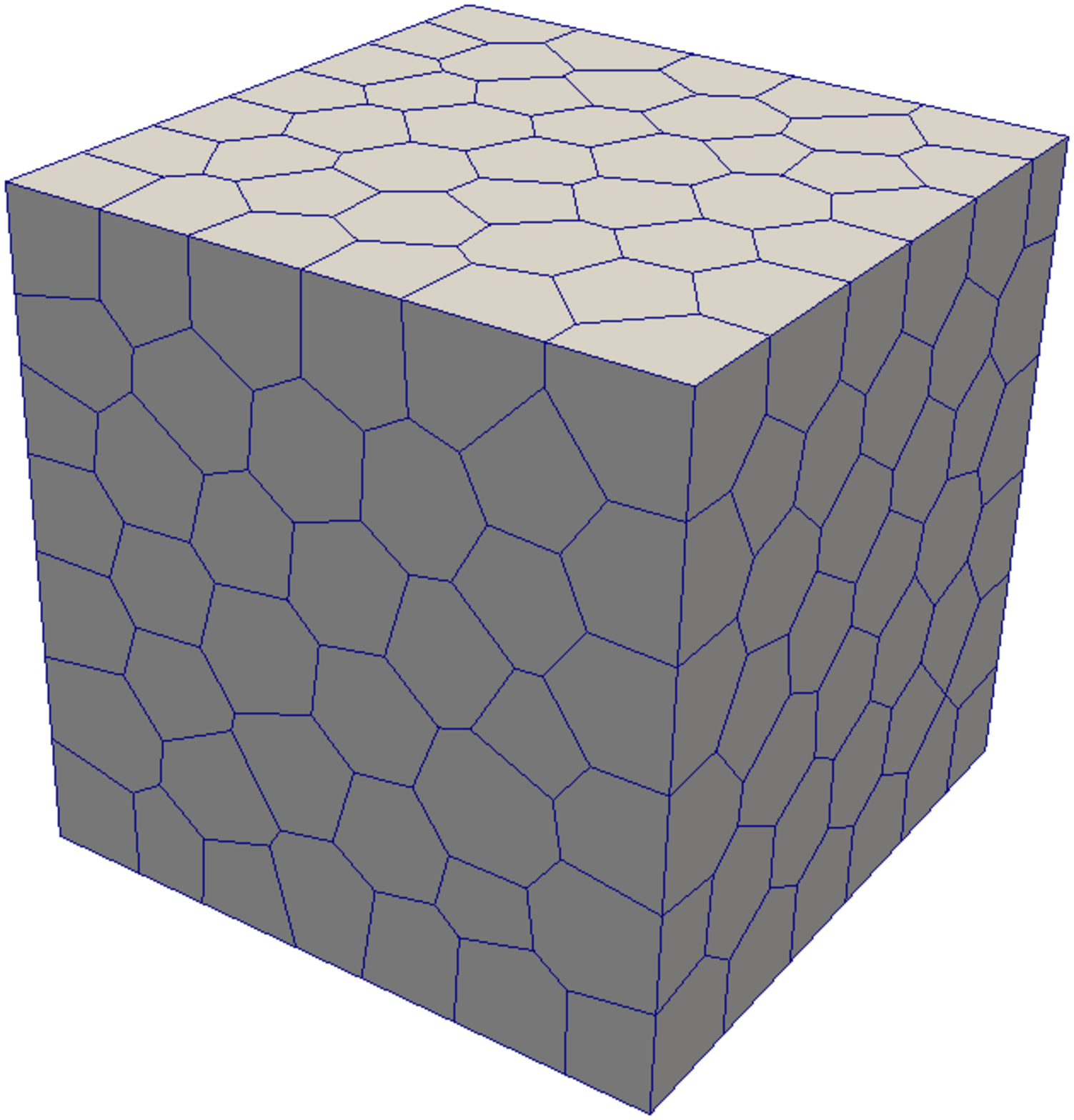}}
\caption{Cube domain discretized with polyhedral elements. Representative meshes containing (a)9, (b) 25, (c) 100 and (d) 300 polyhedra.}
\label{fig:pmesh3d}
\end{figure}

The exact solution to \eref{eqn:problem_strong_form} is $\vm{u} =\hat{\vm{u}}$ in the absence of body forces. The domain is discretized with arbitrary polygonal and polyhedral finite elements. \fref{fig:pmesh2d} and \fref{fig:pmesh3d} shows a few representative meshes used for the two and three dimensional studies, respectively. The errors in the $L^2$ norm and the $H^1$ seminorm for the CS, LS3$n$ schemes and the proposed LS1 one point quadrature are shown in Table \ref{table:linearpatchresults2d} for two-dimensions and in Table \ref{table:linearpatchresults3d} for three dimensions. It can be seen that the proposed one point quadrature scheme passes the linear patch test to machine precision for both polygonal and polyhedral discretizations.

\begin{table}
\renewcommand{\arraystretch}{1.2}
\caption{Error in the $L^2$ norm and $H^1$ seminorm for the two-dimensional linear patch test.}
\centering
\begin{tabular}{lccrccrcc}
\hline
Mesh & \multicolumn{2}{c}{CS} && \multicolumn{2}{c}{LS3$n$-2D} && \multicolumn{2}{c}{LS1-2D}\\
\cline{2-3}\cline{5-6}\cline{8-9}
(c.f. \fref{fig:pmesh2d}) & $L^2$ & $H^1$ && $L^2$ & $H^1$ && $L^2$ & $H^1$ \\
\hline
a &  1.73$\times$10$^{-07}$ & 2.33$\times$10$^{-05}$ && 5.38$\times$10$^{-14}$ & 2.84$\times$10$^{-11}$ && 8.38$\times$10$^{-15}$ & 2.96$\times$10$^{-13}$\\
b & 1.70$\times$10$^{-07}$ & 3.41$\times$10$^{-05}$  && 1.93$\times$10$^{-13}$ & 4.43$\times$10$^{-11}$ && 7.62$\times$10$^{-14}$ & 4.79$\times$10$^{-12}$ \\
c &  7.20$\times$10$^{-07}$ & 2.26$\times$10$^{-04}$ && 2.01$\times$10$^{-13}$ & 7.01$\times$10$^{-11}$ && 1.43$\times$10$^{-13}$ & 1.28$\times$10$^{-11}$ \\
d & 7.42$\times$10$^{-07}$ & 2.58$\times$10$^{-04}$ && 2.96$\times$10$^{-13}$ & 1.02$\times$10$^{-10}$ && 2.71$\times$10$^{-13}$ & 2.76$\times$10$^{-11}$\\
\hline
\end{tabular}
\label{table:linearpatchresults2d}
\end{table}

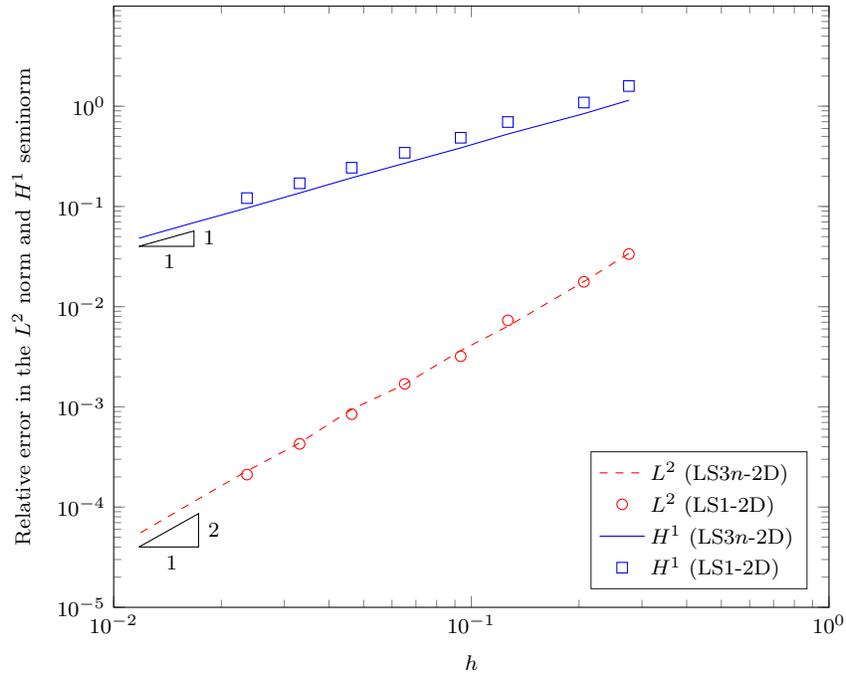
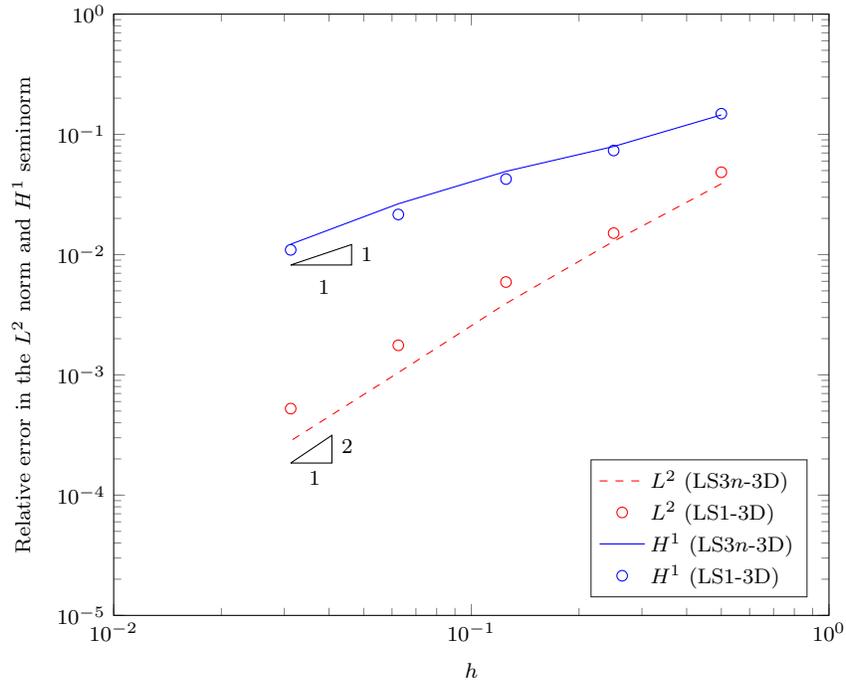
\begin{figure}
\centering
\newlength\figureheight 
\newlength\figurewidth 
\setlength\figureheight{8cm} 
\setlength\figurewidth{10cm}
\subfigure[Two dimensional domain]{
%
%
\begin{tikzpicture}

\begin{axis}[%
width=0.95092\figurewidth,
height=\figureheight,
at={(0\figurewidth,0\figureheight)},
scale only axis,
xmode=log,
xmin=0.01,
xmax=1,
xtick={0.01,  0.1,    1},
xminorticks=true,
xlabel={$h$},
ymode=log,
ymin=1e-05,
ymax=10,
ytick={ 1e-05, 0.0001,  0.001,   0.01,    0.1,      1},
yminorticks=true,
ylabel={Relative error in the $L^2$ norm and $H^1$ seminorm},
legend style={legend cell align=left,align=left,draw=white!15!black},
legend pos = south east
]

\addplot [color=red,dashed]
  table[row sep=crcr]{%
0.275813325481089	0.03398384387782\\
0.206398458721136	0.017769315630541\\
0.126511349842314	0.006375852110579\\
0.0933723503783107	0.003644476687714\\
0.065127276628542	0.001681926480057\\
0.0462983042344264	0.000945417821766958\\
0.0331086783517406	0.000436343249123888\\
0.023584830935152	0.000229194525755016\\
0.0165694147427132	0.000111018629583018\\
0.0117675372744554	5.42965283422432e-05\\
};
\addlegendentry{$L^2$ (LS3$n$-2D)};

\addplot [color=red,only marks,mark=o,mark options={solid}]
  table[row sep=crcr]{%
0.275813325481089	0.0335\\
0.206398458721136	0.0177\\
0.126511349842314	0.0073\\
0.0933723503783107	0.0032\\
0.065127276628542	0.0017\\
0.0462983042344264	0.00084587\\
0.0331086783517406	0.00042765\\
0.023584830935152	0.0002111\\
};
\addlegendentry{$L^2$ (LS1-2D)};

\addplot [color=blue,solid]
  table[row sep=crcr]{%
0.275813325481089	1.14757044668783\\
0.206398458721136	0.843274240727891\\
0.126511349842314	0.52839888314269\\
0.0933723503783107	0.384069226737528\\
0.065127276628542	0.269758713868738\\
0.0462983042344264	0.19309035802328\\
0.0331086783517406	0.136291797527042\\
0.023584830935152	0.096540951598009\\
0.0165694147427132	0.068202755819459\\
0.0117675372744554	0.048244263738452\\
};
\addlegendentry{$H^1$ (LS3$n$-2D)};

\addplot [color=blue,only marks,mark=square,mark options={solid}]
  table[row sep=crcr]{%
0.275813325481089	1.591\\
0.206398458721136	1.0884\\
0.126511349842314	0.6958\\
0.0933723503783107	0.4848\\
0.065127276628542	0.3434\\
0.0462983042344264	0.2436\\
0.0331086783517406	0.1703\\
0.023584830935152	0.1211\\
};
\addlegendentry{$H^1$ (LS1-2D)};

\addplot [color=black,solid,forget plot]
  table[row sep=crcr]{%
0.0117675372744554	0.04\\
0.0167675372744554	0.04\\
0.0167675372744554	0.0569959096228363\\
0.0117675372744554	0.04\\
};
\node[right, align=left, inner sep=0mm, text=black]
at (axis cs:0.0137675372744554,0.028,0) {$1$};
\node[right, align=left, inner sep=0mm, text=black]
at (axis cs:0.0177735895109227,0.0487686014804363,0) {$1$};
\addplot [color=black,solid,forget plot]
  table[row sep=crcr]{%
0.0117675372744554	4e-05\\
0.0172675372744554	4e-05\\
0.0172675372744554	8.61290447234452e-05\\
0.0117675372744554	4e-05\\
};
\node[right, align=left, inner sep=0mm, text=black]
at (axis cs:0.0139675372744554,2.8e-05,0) {$1$};
\node[right, align=left, inner sep=0mm, text=black]
at (axis cs:0.0183035895109227,5.91397877163212e-05,0) {$2$};
\end{axis}
\end{tikzpicture}%
}
\subfigure[Three dimensional domain]{
%
%
\begin{tikzpicture}

\begin{axis}[%
width=0.95092\figurewidth,
height=\figureheight,
at={(0\figurewidth,0\figureheight)},
scale only axis,
xmode=log,
xmin=0.01,
xmax=1,
xtick={0.01,  0.1,    1},
xminorticks=true,
xlabel={$h$},
ymode=log,
ymin=1e-05,
ymax=1,
ytick={ 1e-05, 0.0001,  0.001,   0.01,    0.1,      1},
yminorticks=true,
ylabel={Relative error in the $L^2$ norm and $H^1$ seminorm},
legend style={at={(0.97,0.03)},anchor=south east,legend cell align=left,align=left,draw=white!15!black}
]
\addplot [color=red,dashed]
  table[row sep=crcr]{%
0.5	0.038821568082792\\
0.25	0.012906053738006\\
0.125	0.003923804883409\\
0.0625	0.0010449469075704\\
0.03125	0.000282669452814865\\
};
\addlegendentry{$L^2$ (LS3$n$-3D)};

\addplot [color=red,only marks,mark=o,mark options={solid}]
  table[row sep=crcr]{%
0.5	0.048410475288474\\
0.25	0.015094999796221\\
0.125	0.005919238362217\\
0.0625	0.001761737506258\\
0.03125	0.000524344324561663\\
};
\addlegendentry{$L^2$ (LS1-3D)};

\addplot [color=blue,solid]
  table[row sep=crcr]{%
0.5	0.14514106472845\\
0.25	0.079306725432887\\
0.125	0.049247194183739\\
0.0625	0.026463905630985\\
0.03125	0.0122087720638\\
};
\addlegendentry{$H^1$ (LS3$n$-3D)};

\addplot [color=blue,only marks,mark=o,mark options={solid}]
  table[row sep=crcr]{%
0.5	0.148517069183995\\
0.25	0.073406696073821\\
0.125	0.042498915500584\\
0.0625	0.021576418942673\\
0.03125	0.01095420551104\\
};
\addlegendentry{$H^1$ (LS1-3D)};

\addplot [color=black,solid,forget plot]
  table[row sep=crcr]{%
0.03125	0.000185\\
0.04075	0.000185\\
0.04075	0.00031457696\\
0.03125	0.000185\\
};
\node[right, align=left, inner sep=0mm, text=black]
at (axis cs:0.03505,0.00013875,0) {$1$};
\node[right, align=left, inner sep=0mm, text=black]
at (axis cs:0.043195,0.000254905991111111,0) {$2$};
\addplot [color=black,solid,forget plot]
  table[row sep=crcr]{%
0.03125	0.0082\\
0.04625	0.0082\\
0.04625	0.012136\\
0.03125	0.0082\\
};
\node[right, align=left, inner sep=0mm, text=black]
at (axis cs:0.03725,0.00533,0) {$1$};
\node[right, align=left, inner sep=0mm, text=black]
at (axis cs:0.049025,0.0100670769230769,0) {$1$};
\end{axis}
\end{tikzpicture}%
}
\caption{Convergence results for the quadratic patch test. The domain is discretized with arbitrary polytopes. The new integration scheme delivers optimal convergence rates in both the $L^2$ norm and $H^1$ seminorm with three times as many integration points per element as the standard approach.}
\label{fig:quadpatchlinearConveResults}
\end{figure}

\begin{table}
\renewcommand{\arraystretch}{1.2}
\caption{Error in the $L^2$ norm and $H^1$ seminorm for the three-dimensional linear patch test.}
\centering
\begin{tabular}{lccrcc}
\hline
Mesh & \multicolumn{2}{c}{LS3$n$-3D} && \multicolumn{2}{c}{LS1-3D} \\
\cline{2-3}\cline{5-6}
(c.f. \fref{fig:pmesh3d}) & $L^2$ & $H^1$ && $L^2$ & $H^1$ \\
\hline
$a$ & 2.03$\times$10$^{-12}$ & 3.34$\times$10$^{-10}$ && 2.98$\times$10$^{-11}$ & 2.23$\times$10$^{-10}$\\
$b$ & 1.92$\times$10$^{-12}$ & 1.75$\times$10$^{-10}$ && 7.38$\times$10$^{-10}$ & 5.56$\times$10$^{-09}$\\
$c$ & 2.66$\times$10$^{-12}$ & 4.93$\times$10$^{-10}$&& 2.08$\times$10$^{-10}$ & 2.13$\times$10$^{-09}$ \\
$d$ & 3.21$\times$10$^{-12}$ & 3.11$\times$10$^{-10}$ && 7.73$\times$10$^{-10}$ & 1.28$\times$10$^{-09}$  \\
\hline
\end{tabular}
\label{table:linearpatchresults3d}
\end{table}
Next, to study the convergence properties of the proposed technique, the following higher order displacements are prescribed on the boundaries: 
\begin{equation}
\begin{pmatrix} \hat{u} \\ \hat{v} \end{pmatrix} = \begin{pmatrix} 0.1x^2+0.1xy+0.2y^2 \\ 0.05x^2+0.15xy + 0.1y^2 \end{pmatrix},
\end{equation}
in the two-dimensional case and the following in the three-dimensional case:
\begin{equation}
\begin{pmatrix} \hat{u} \\ \hat{v} \\ \hat{w} \end{pmatrix} = \begin{pmatrix} 0.1+0.2x+0.2x+0.1z+0.15x^2+0.2y^2+0.1z^2+0.15xy+0.1yz+0.1zx \\ 0.15+0.1x+0.1y+0.2z+0.2x^2+0.15y^2+0.1z^2+0.2xy+0.1yz+0.2zx\\ 0.15+0.15x+0.2y+0.1z+0.15x^2+0.1y^2+0.2z^2+0.1xy+0.2yz+0.15zx \end{pmatrix}
\end{equation}
The exact solution to \eref{eqn:problem_strong_form} is $\vm{u} = \hat{\vm{u}} $ when the body is subjected to the body forces:
\begin{equation}
\mathbf{b} = \begin{pmatrix} -0.2\mat{C}(1,1) -0.15\mat{C}(1,2)-0.55\mat{C}(3,3) \\ -0.1\mat{C}(1,2)-0.2\mat{C}(2,2)-0.2\mat{C}(3,3) \end{pmatrix},
\end{equation}
in two-dimensions and
\begin{equation}
\mathbf{b} = \begin{pmatrix} -0.3\mat{C}(1,1) -0.2\mat{C}(1,2)-0.15\mat{C}(1,3) -0.6\mat{C}(4,4) - 0.35\mat{C}(6,6) \\ -0.15\mat{C}(1,2)-0.3\mat{C}(2,2)-0.2\mat{C}(2,3) - 0.55\mat{C}(4,4) - 0.4\mat{C}(5,5) \\ 0.1\mat{C}(1,3) - 0.1\mat{C}(2,3) - 0.4\mat{C}(3,3) - 0.3\mat{C}(5,5) - 0.4\mat{C}(6,6) \end{pmatrix}
\end{equation}
in three dimensions, where $\mat{C}$ is the constitutive matrix. \fref{fig:quadpatchlinearConveResults} shows the convergence rates when the domain is discretized with polyhedral linear elements. It can be inferred that the proposed one point quadrature scheme yields optimal convergence rates.

\subsection{Thick cantilever beam under end shear}
\begin{figure}[htpb]
\centering
\includegraphics[scale=1]{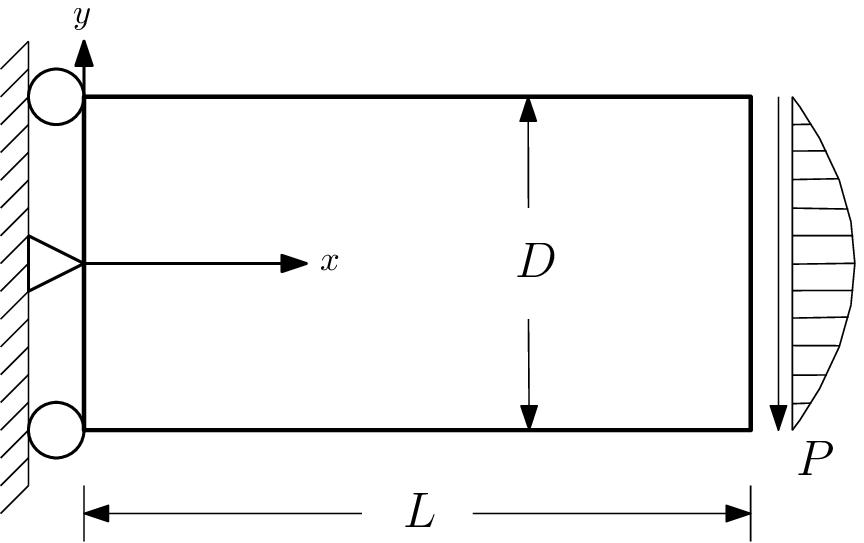}
\caption{Geometry and boundary conditions for the two dimensional cantilever beam problem.}
\label{fig:cantileverfig}
\end{figure}

\begin{figure}[htpb]
\includegraphics[width=1\textwidth]{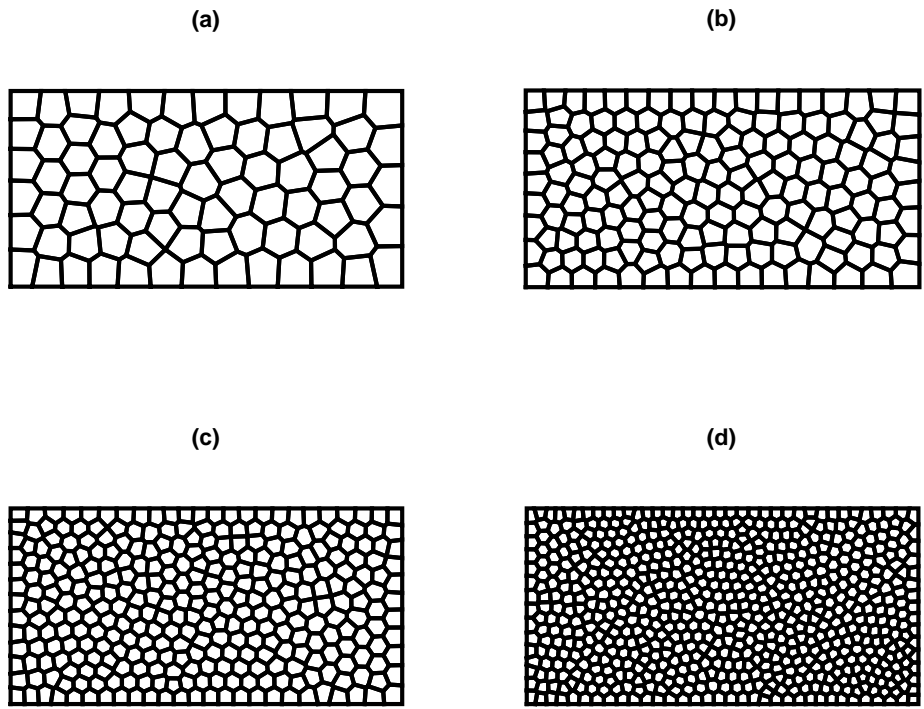}
\caption{Sample meshes for the two dimensional cantilever beam problem containing: (a) 80, (b) 160, (c) 320 and (d) 640 polygons.}
\label{fig:cantmesh}
\end{figure}

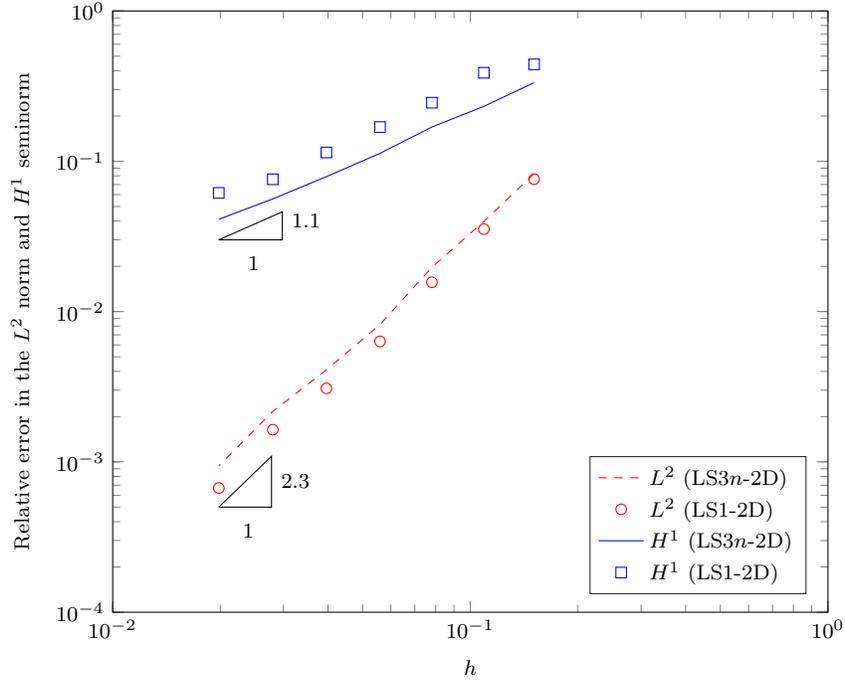
\begin{figure}
\centering
\setlength\figureheight{8cm} 
\setlength\figurewidth{10cm}
%
%
\begin{tikzpicture}

\begin{axis}[%
width=0.95092\figurewidth,
height=\figureheight,
at={(0\figurewidth,0\figureheight)},
scale only axis,
xmode=log,
xmin=0.01,
xmax=1,
xminorticks=true,
xlabel={$h$},
ymode=log,
ymin=0.0001,
ymax=1,
yminorticks=true,
ylabel={Relative error in the $L^2$ norm and $H^1$ seminorm},
legend style={legend cell align=left,align=left,draw=white!15!black},
legend pos = south east
]

\addplot [color=red,dashed]
  table[row sep=crcr]{%
0.150755672288882	0.082584606755026\\
0.109108945117996	0.039771977105347\\
0.078086880944303	0.01971605464504\\
0.0559016994374947	0.008192631465247\\
0.0395903791232448	0.004099044019674\\
0.0280386077046022	0.002168791404893\\
0.0198029508595335	0.000940236645391731\\
};
\addlegendentry{$L^2$ (LS3$n$-2D)};

\addplot [color=red,only marks,mark=o,mark options={solid}]
  table[row sep=crcr]{%
0.150755672288882	0.076017674122064\\
0.109108945117996	0.035369060654454\\
0.078086880944303	0.015661100905351\\
0.0559016994374947	0.006331204929675\\
0.0395903791232448	0.003085285369967\\
0.0280386077046022	0.001637626838452\\
0.0198029508595335	0.000670244907293391\\
};
\addlegendentry{$L^2$ (LS1-2D)};

\addplot [color=blue,solid]
  table[row sep=crcr]{%
0.150755672288882	0.333374254304147\\
0.109108945117996	0.231856134064404\\
0.078086880944303	0.168965632100084\\
0.0559016994374947	0.112864449467272\\
0.0395903791232448	0.079008210331788\\
0.0280386077046022	0.056187065898335\\
0.0198029508595335	0.041086635887043\\
};
\addlegendentry{$H^1$ (LS3$n$-2D)};

\addplot [color=blue,only marks,mark=square,mark options={solid}]
  table[row sep=crcr]{%
0.150755672288882	0.441149315954171\\
0.109108945117996	0.387676066386607\\
0.078086880944303	0.245177063533072\\
0.0559016994374947	0.168855147821359\\
0.0395903791232448	0.114389173094795\\
0.0280386077046022	0.075860212187247\\
0.0198029508595335	0.06158184168984\\
};
\addlegendentry{$H^1$ (LS1-2D)};

\addplot [color=black,solid,forget plot]
  table[row sep=crcr]{%
0.0198029508595335	0.0005\\
0.0278029508595335	0.0005\\
0.0278029508595335	0.00109118971445855\\
0.0198029508595335	0.0005\\
};
\node[right, align=left, inner sep=0mm, text=black]
at (axis cs:0.0230029508595335,0.00035,0) {$1$};
\node[right, align=left, inner sep=0mm, text=black]
at (axis cs:0.0294711279111055,0.000742486603213012,0) {$2.3$};
\addplot [color=black,solid,forget plot]
  table[row sep=crcr]{%
0.0198029508595335	0.03\\
0.0298029508595335	0.03\\
0.0298029508595335	0.0460815500951261\\
0.0198029508595335	0.03\\
};
\node[right, align=left, inner sep=0mm, text=black]
at (axis cs:0.0238029508595335,0.021,0) {$1$};
\node[right, align=left, inner sep=0mm, text=black]
at (axis cs:0.0315911279111055,0.0402403444655836,0) {$1.1$};
\end{axis}
\end{tikzpicture}%
\caption{Convergence of the relative error in the $L^2$ norm and the $H^1$ seminorm with mesh refinement for a two-dimensional cantilever beam subjected to end shear. It is inferred that the proposed integration scheme yields optimal convergence rates.}
\label{fig:cantiConveResults_L2H1}
\end{figure}

In this example, a two-dimensional cantilever beam subjected to a parabolic shear load at the free end is examined, as shown in \fref{fig:cantileverfig}. The geometry of the cantilever is $L=$ 10 m and  $D=$ 2 m. The material properties are: Young's modulus, $E=$ 3$\times 10^7$ N/m$^2$, Poisson's ratio $\nu=$ 0.25 and the parabolic shear force is $P=$ 150 N. The exact solution for the displacement field is given by~\cite{Timoshenko1970}:
\begin{align}
u(x,y) &= \frac{P y}{6 EI} \left[ (9L-3x)x + (2+\nu) \left( y^2 - \frac{D^2}{4} \right) \right], \nonumber \\
v(x,y) &= -\frac{P}{6 EI} \left[ 3 \nu y^2(L-x) + (4+5\nu) \frac{D^2x}{4} + (3L-x)x^2 \right].
\label{eqn:cantisolution}
\end{align}
where $I = D^3/12$ is the second area moment. A state of plane stress is considered. \fref{fig:cantmesh} shows few sample polygonal meshes. The numerical convergence of the relative error in the $L^2$ norm and the $H^1$ seminorm is shown in \fref{fig:cantiConveResults_L2H1}. It can be seen that the proposed one point integration rule yields optimal convergence rate in both the $L^2$ norm and the $H^1$ seminorm. With mesh refinement the solution approaches the analytical solution asymptotically. It is further noted that the proposed integration rule yields similar results when compared to the recently proposed integration rule~\cite{francisa.ortiz-bernardin2017} that employs 3$n$ integration point per element (see \fref{fig:subcellRepre}a).

\subsection{Three dimensional cantilever beam under end torsion}
\begin{figure}[htpb]
\centering
\includegraphics[scale=1]{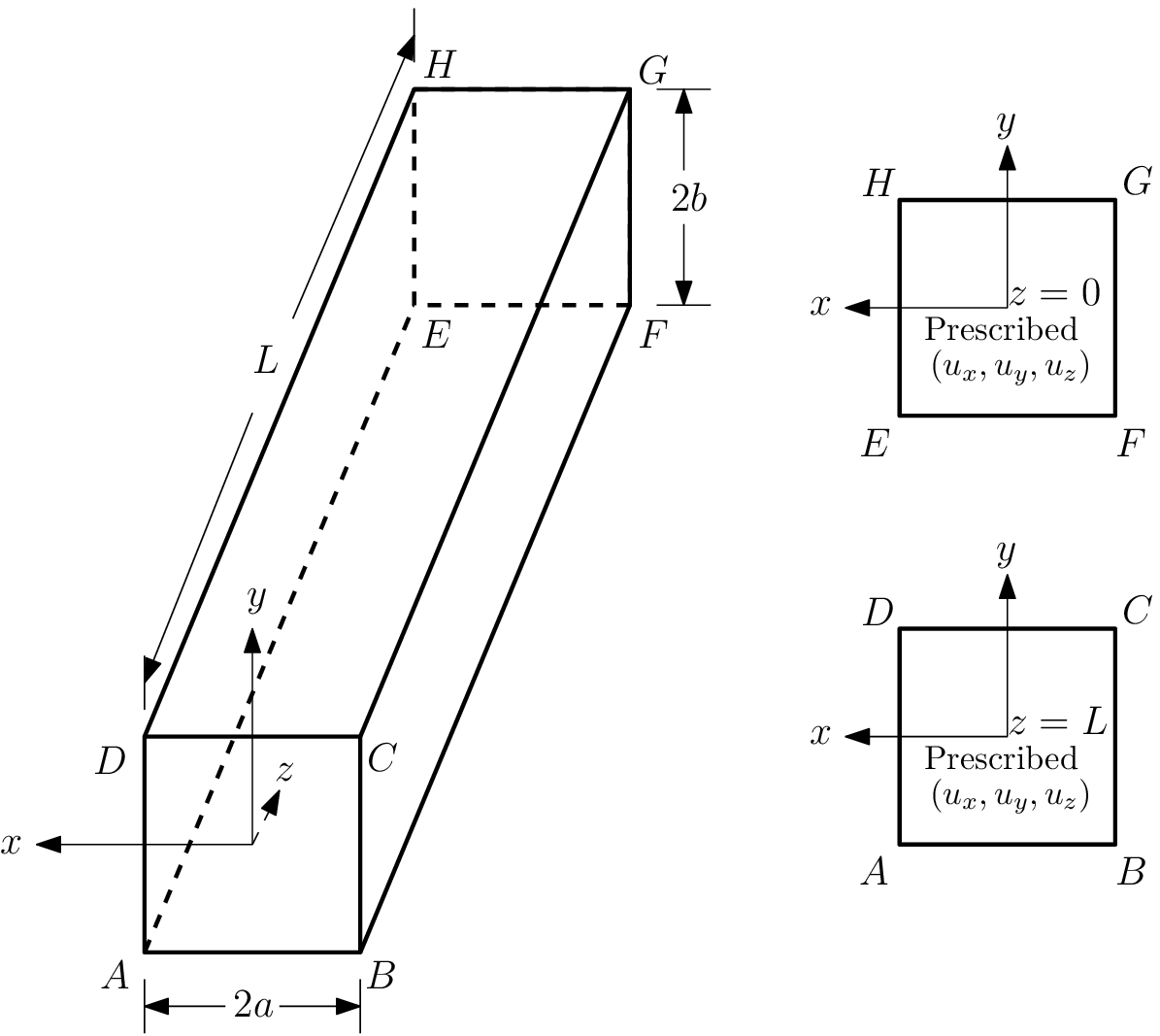}
\caption{Cantilever beam: (a) Geometry, length $L$ and rectangular cross-section of width $2a$ and height $2b$. For the present study, the following dimensions are considered: $L=$ 5, $a=b=$ 1.}
\label{fig:canti3dbeam}
\end{figure}

\begin{figure}
\centering
\subfigure[]{\includegraphics[scale=0.17]{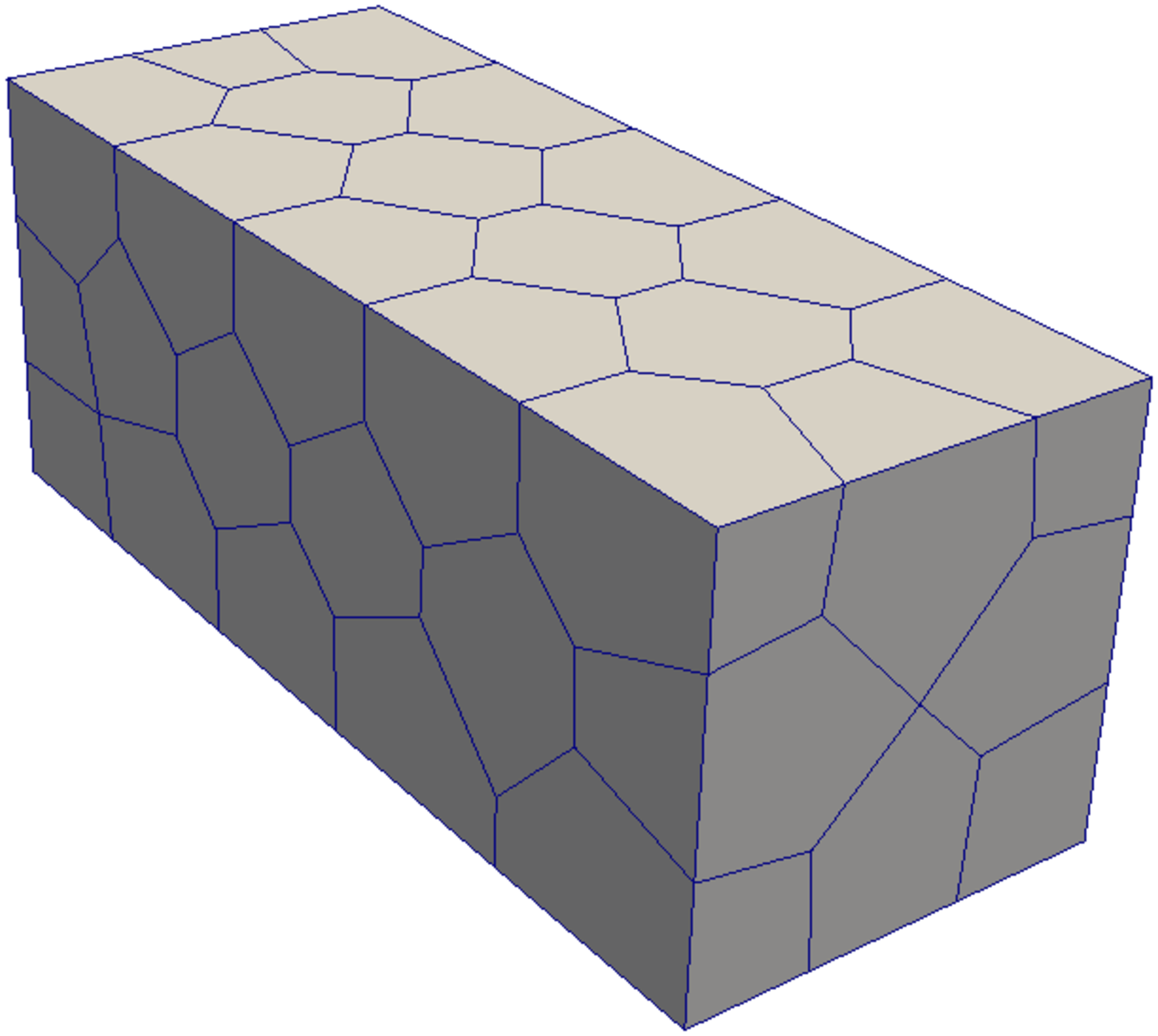}}
\subfigure[]{\includegraphics[scale=0.17]{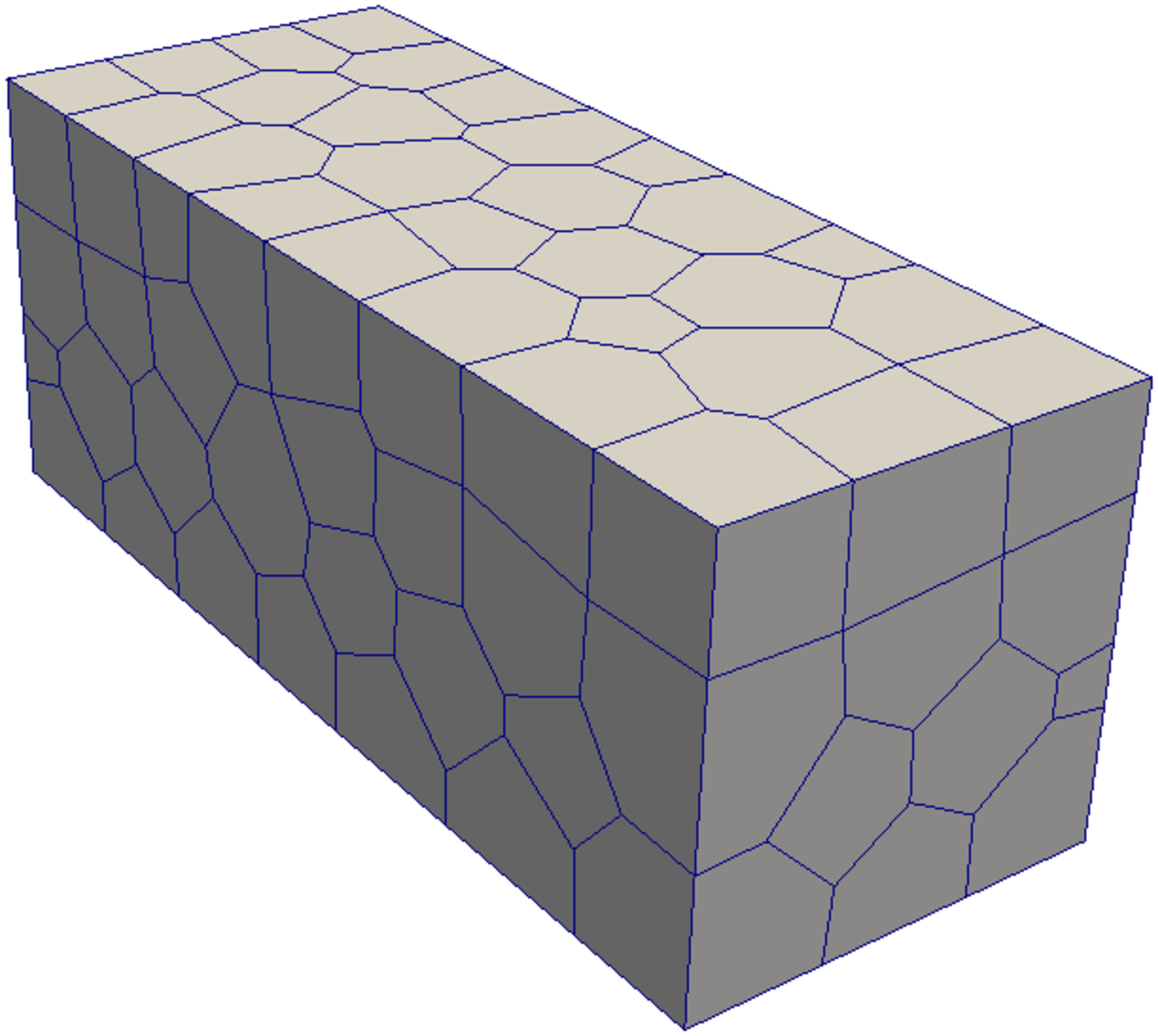}}
\subfigure[]{\includegraphics[scale=0.17]{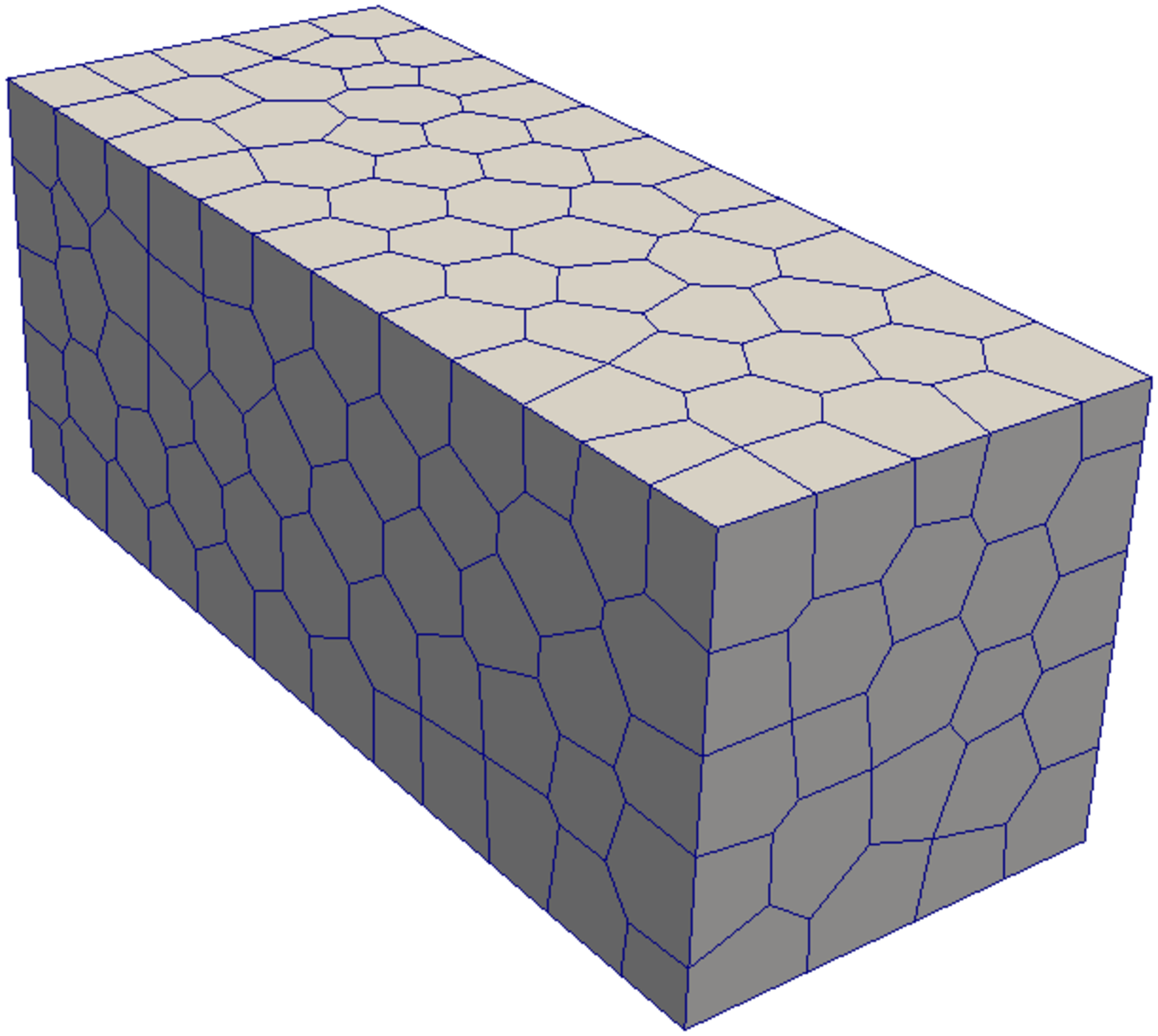}}
\subfigure[]{\includegraphics[scale=0.17]{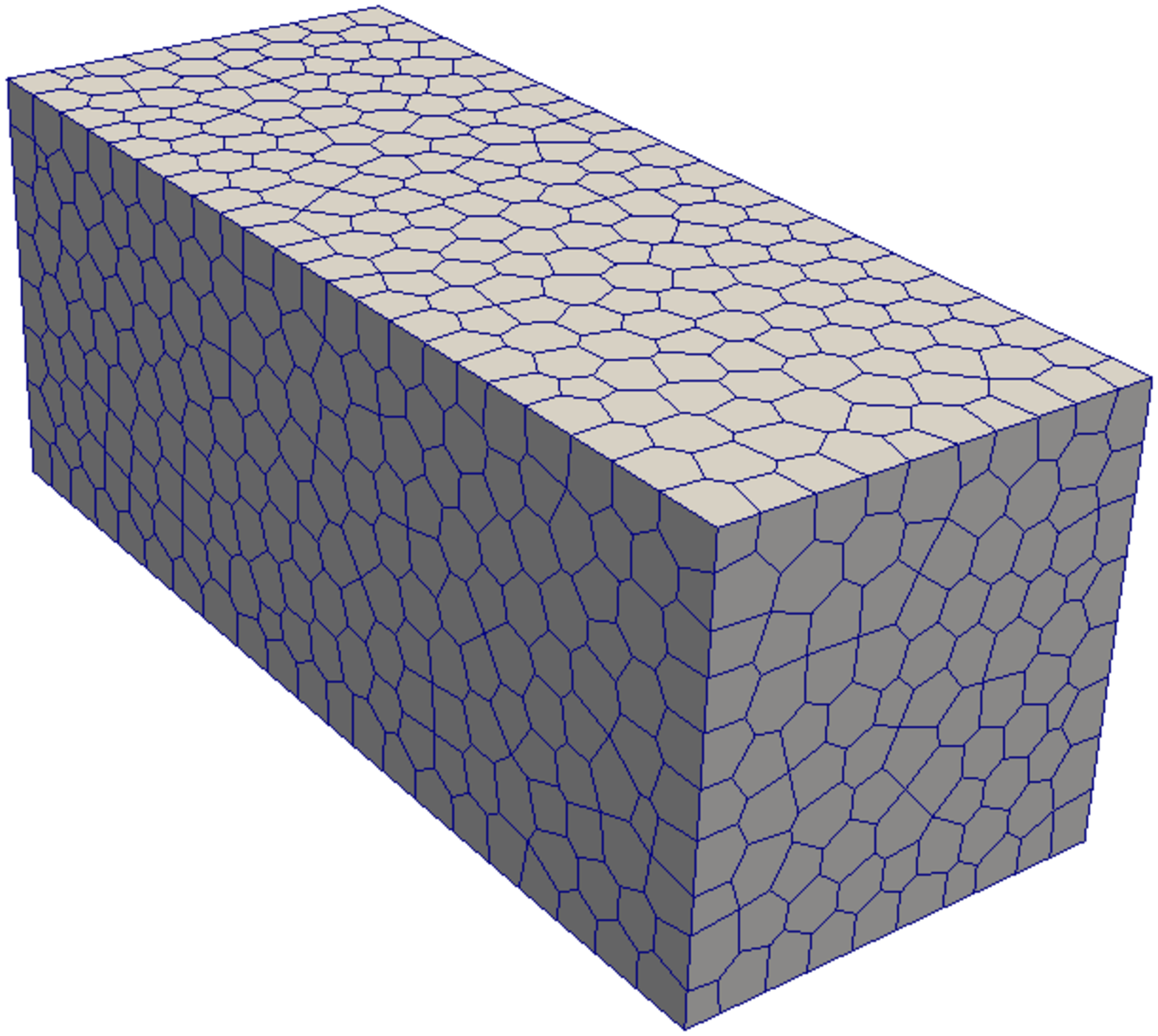}}
\caption{Sample meshes for the three dimensional cantilever beam problem containing (a) 50,  (b) 100, (c) 300 and (d) 2000 polyhedra.}
\label{fig:cantirandmesh}
\end{figure}

\begin{figure}[htpb]
\centering
\setlength\figureheight{8cm} 
\setlength\figurewidth{10cm}
%
%
\begin{tikzpicture}

\begin{axis}[%
width=0.95092\figurewidth,
height=\figureheight,
at={(0\figurewidth,0\figureheight)},
scale only axis,
xmode=log,
xmin=0.01,
xmax=1,
xtick={0.01,  0.1,    1},
xminorticks=true,
xlabel={Maximum edge size $$h$$},
ymode=log,
ymin=0.001,
ymax=1,
ytick={ 1e-05, 0.0001,  0.001,   0.01,    0.1,      1},
yminorticks=true,
ylabel={Relative error in the $L^2$ norm and $H^1$ seminorm},
legend style={at={(0.97,0.03)},anchor=south east,legend cell align=left,align=left,draw=white!15!black}
]
\addplot [color=red,dashed]
  table[row sep=crcr]{%
0.110959112137338	0.110757285379062\\
0.0867266587617825	0.070196368874268\\
0.0682392773497722	0.04099518640725\\
0.0537379988149689	0.025193972510603\\
0.0423195629003272	0.015022655414559\\
0.0334064928071621	0.009055170978243\\
};
\addlegendentry{$L^2$ (LS3$n$-3D)};

\addplot [color=red,only marks,mark=o,mark options={solid}]
  table[row sep=crcr]{%
0.110959112137338	0.113563260652218\\
0.0867266587617825	0.069625998963438\\
0.0682392773497722	0.039958349954481\\
0.0537379988149689	0.023547884447573\\
0.0423195629003272	0.014361641411575\\
0.0334064928071621	0.009099740681235\\
};
\addlegendentry{$L^2$ (LS1-3D)};

\addplot [color=blue,solid]
  table[row sep=crcr]{%
0.110959112137338	0.676131968546381\\
0.0867266587617825	0.531644272113483\\
0.0682392773497722	0.412070960785253\\
0.0537379988149689	0.322605271037476\\
0.0423195629003272	0.245960083081151\\
0.0334064928071621	0.190645185511065\\
};
\addlegendentry{$H^1$ (LS3$n$-3D)};

\addplot [color=blue,only marks,mark=square,mark options={solid}]
  table[row sep=crcr]{%
0.110959112137338	0.701592289068987\\
0.0867266587617825	0.518588725072203\\
0.0682392773497722	0.385632227012935\\
0.0537379988149689	0.322557072830983\\
0.0423195629003272	0.252908534975217\\
0.0334064928071621	0.201126794035026\\
};
\addlegendentry{$H^1$ (LS1-3D)};

\addplot [color=black,solid,forget plot]
  table[row sep=crcr]{%
0.0334064928071621	0.00725\\
0.0429064928071621	0.00725\\
0.0429064928071621	0.0119597547178919\\
0.0334064928071621	0.00725\\
};
\node[right, align=left, inner sep=0mm, text=black]
at (axis cs:0.0372064928071621,0.0054375,0) {$1$};
\node[right, align=left, inner sep=0mm, text=black]
at (axis cs:0.0454808823755918,0.00990772327064263,0) {$2$};
\addplot [color=black,solid,forget plot]
  table[row sep=crcr]{%
0.0334064928071621	0.15\\
0.0484064928071621	0.15\\
0.0484064928071621	0.217352176506168\\
0.0334064928071621	0.15\\
};
\node[right, align=left, inner sep=0mm, text=black]
at (axis cs:0.0394064928071621,0.1125,0) {$1$};
\node[right, align=left, inner sep=0mm, text=black]
at (axis cs:0.0513108823755918,0.183438796385564,0) {$1$};
\end{axis}
\end{tikzpicture}%
\caption{Convergence of the relative error in the $L^2$ norm and the $H^1$ seminorm with mesh refinement for the three-dimensional cantilever beam problem subjected to end torsion. It can be seen that the proposed integration rule yields similar results when compared to linear smoothing scheme. The rate of convergence is also optimal in both the $L^2$ norm and in the $H^1$ seminorm.}
\label{fig:cbeam3dL2_torsion}
\end{figure}
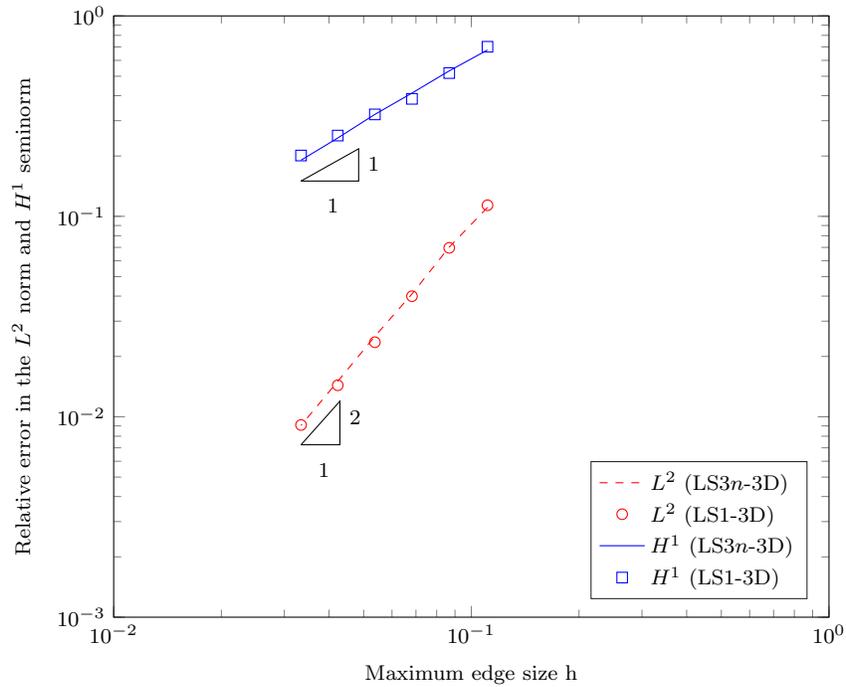

Consider a prismatic cantilever beam with $\Omega: [-1,1] \times [-1,1] \times [0,L]$ (see \fref{fig:canti3dbeam}
(a) for geometry of the domain) subjected to end torsion. The material is assumed to be homogeneous and isotropic with Youngs' modulus, $E=$ 1 N/m$^2$, Poisson's ratio $\nu=$ 0.3 and shear modulus $G=E/(2(1+\nu))$. Two different loading conditions, viz., end shear load and end torsion, are considered here for which analytical solutions are available in the literature. The accuracy and the convergence properties are studied for random closed-pack Voronoi mesh. \fref{fig:cantirandmesh} shows a few representative random Voronoi meshes employed for this study.
The exact displacement solution for this boundary value problem is~\cite{barber2010}:
\begin{align}
u_x &= -\beta yz \nonumber \\
u_y &= \beta xz \nonumber \\
u_z &= \beta \left[ xy + \sum\limits_{n=1}^\infty \frac{32a^2(-1)^n}{\pi^3 (2n-1)^3} \sin \left( (2n-1) \frac{\pi x}{2a} \right) \frac{ \sinh ((2n-1)\frac{\pi y}{2a})}{ \cosh ((2n-1)\frac{\pi y}{2a})} \right]
\label{eqn:tordisp}
\end{align}
where the constant $\beta$ is proportional to the total torque applied to the beam. The exact Cauchy stress field is given by:
\begin{align}
\sigma_{xx} &= \sigma_{xy} = \sigma_{yy} = \sigma_{zz} = 0 \nonumber \\
\sigma_{xz} &= G \beta \sum\limits_{n=1}^\infty \frac{ 16a(-1)^n}{\pi^2 (2n-1)^2} \cos \left( (2n-1) \frac{\pi x}{2a} \right) \frac{ \sinh ((2n-1)\frac{\pi y}{2a})}{ \cosh ((2n-1)\frac{\pi y}{2a})} \nonumber \\
\sigma_{yz} &= G\beta \left[ 2x + \sum\limits_{n=1}^\infty \frac{ 16a(-1)^n}{\pi^2 (2n-1)^2} \sin \left( (2n-1) \frac{\pi x}{2a} \right) \frac{ \cosh ((2n-1)\frac{\pi y}{2a})}{ \cosh ((2n-1)\frac{\pi y}{2a})}  \right]
\label{eqn:torstress}
\end{align}
The infinite series in \erefs{eqn:tordisp} - (\ref{eqn:torstress}) is truncated at $n=$ 40. The exact solution for the displacement is prescribed on the surface at $z=L$ and at $z=0$, surface tractions are applied at the rest of the boundary, which are consistent with the exact stress field. The convergence of the proposed technique over arbitrary polyhedron with mesh refinement is studied. The error in the $L^2$ and the $H^1$ seminorm is shown in \fref{fig:cbeam3dL2_torsion} and it can be seen that the proposed approach yields optimal convergence rates. The results from the present approach is compared with the linear smoothing technique that employs 4 integration points per tetrahedron.

\subsection{Three dimensional L-shaped block}
\begin{figure}[htpb]
\centering
\subfigure[A 3D block problem: Geometry and boundary conditions.]{\includegraphics[width=0.45\textwidth]{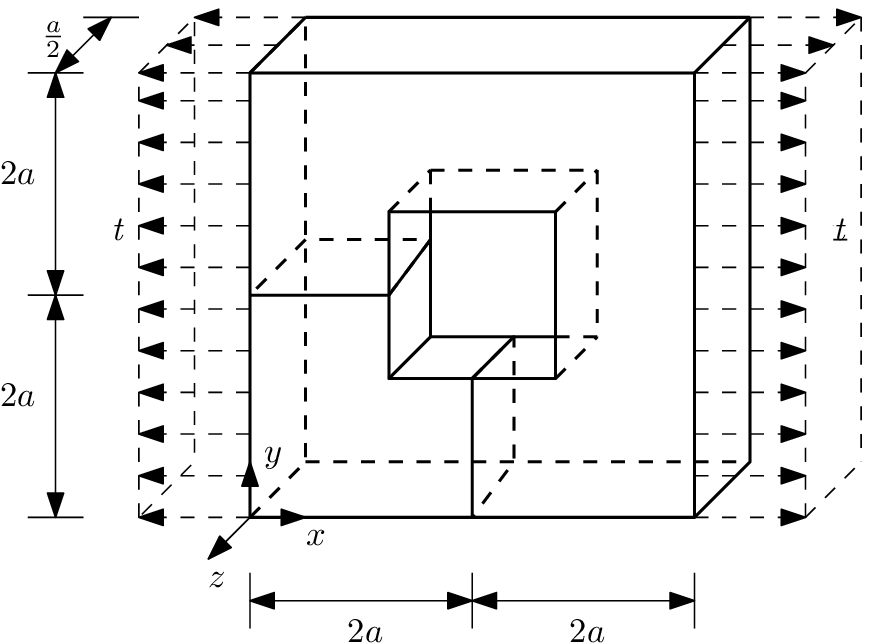}}
\hspace{1cm}
\subfigure[An L-shaped quarter model: Geometry and boundary conditions.]{\includegraphics[width=0.4\textwidth]{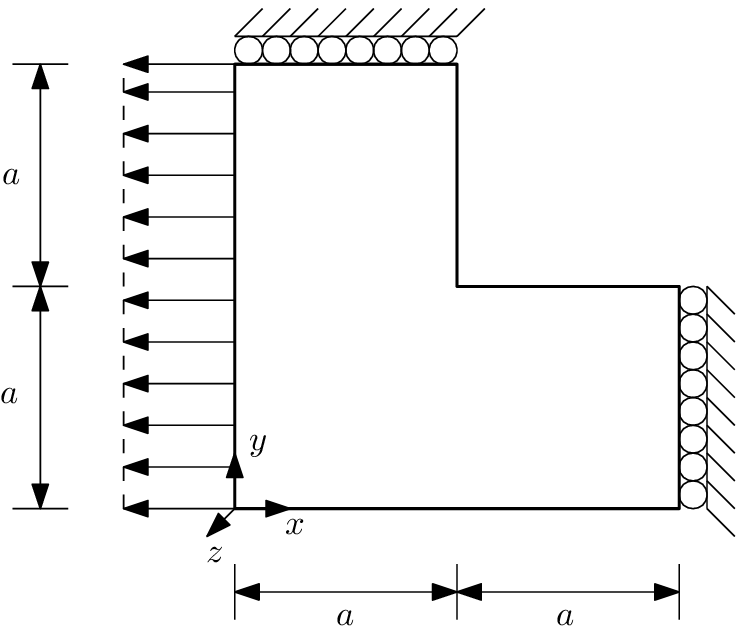}}
\caption{Three dimensional block and an L-shaped quarter model.}
\label{fig:Lshape}
\end{figure}

\begin{figure}[htpb]
\centering
\subfigure[]{\includegraphics[width=0.47\textwidth]{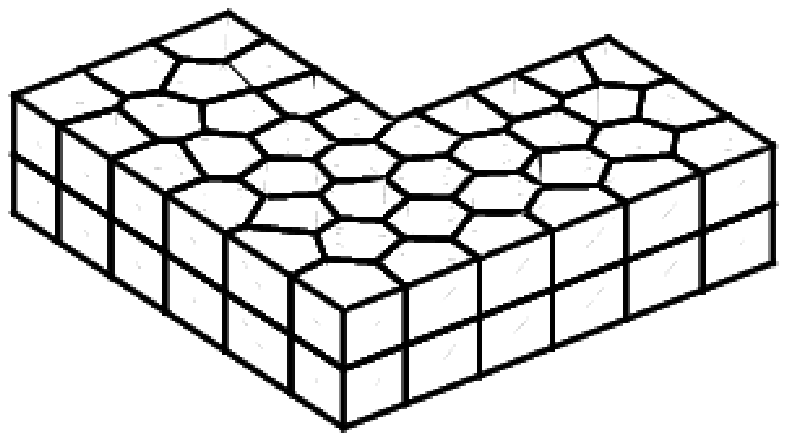}}
\subfigure[]{\includegraphics[width=0.47\textwidth]{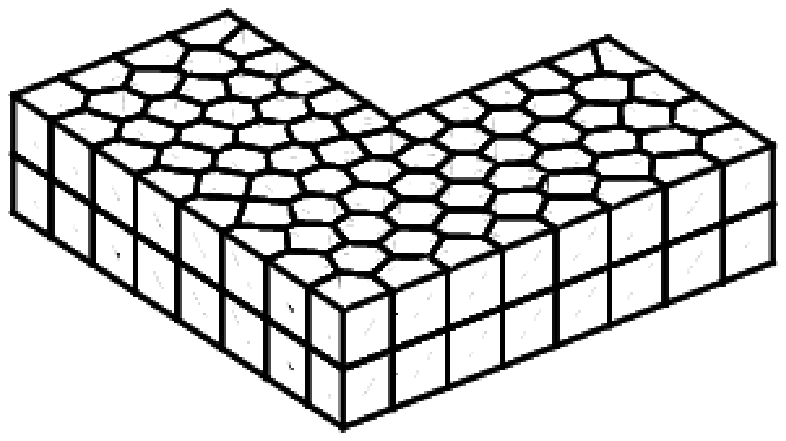}}
\subfigure[]{\includegraphics[width=0.47\textwidth]{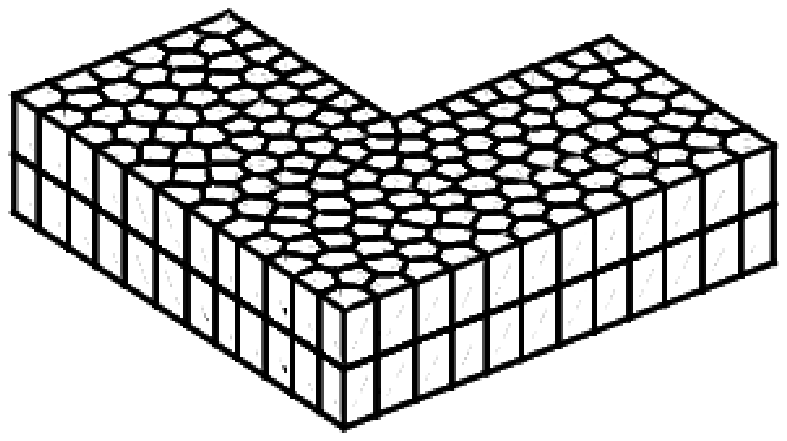}}
\subfigure[]{\includegraphics[width=0.47\textwidth]{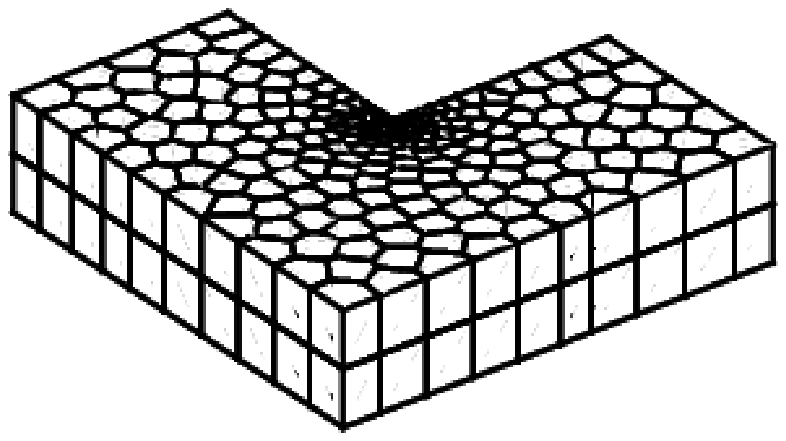}}
\caption{Few sample meshes of L-shaped quarter model containing: a) 40 elements b) 80 elements c)160 elements and d) 320 elements.}
\label{fig:3dmesh}
\end{figure}

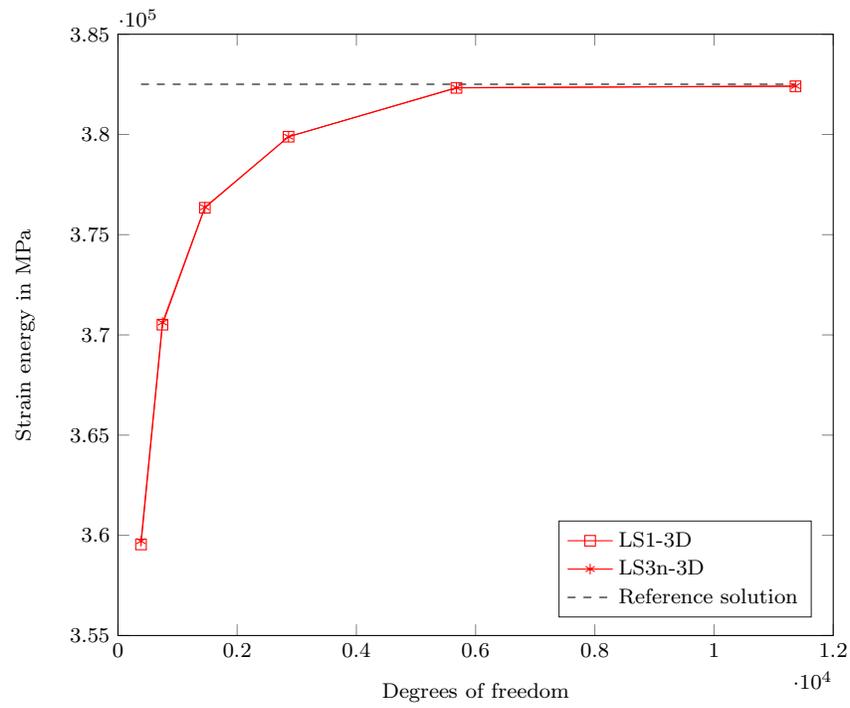
\begin{figure}[htpb]
\centering
\setlength\figureheight{8cm} 
\setlength\figurewidth{10cm}
%
%
\begin{tikzpicture}

\begin{axis}[%
width=0.95092\figurewidth,
height=\figureheight,
at={(0\figurewidth,0\figureheight)},
scale only axis,
xmin=0,
xmax=12000,
xlabel={Degrees of freedom},
ymin=355000,
ymax=385000,
ylabel={Strain energy in MPa},
legend style={at={(0.97,0.03)},anchor=south east,legend cell align=left,align=left,draw=white!15!black}
]
\addplot [color=red,solid,mark=square,mark options={solid}]
  table[row sep=crcr]{%
387	359543.950980508\\
747	370511.038069552\\
1458	376341.671922878\\
2862	379887.605421102\\
5679	382323.134668055\\
11367	382403.695747692\\
};
\addlegendentry{LS1-3D};

\addplot [color=red,solid,mark=asterisk,mark options={solid}]
  table[row sep=crcr]{%
387	359709.372951002\\
747	370616.877657546\\
1458	376370.688840161\\
2862	379887.782227358\\
5679	382335.561215619\\
11367	382414.855103007\\
};
\addlegendentry{LS3n-3D};

\addplot [color=black,dashed]
  table[row sep=crcr]{%
387	382505\\
747	382505\\
1458	382505\\
2862	382505\\
5679	382505\\
11367	382505\\
};
\addlegendentry{Reference solution};

\end{axis}
\end{tikzpicture}%
\caption{Convergence of the strain energy of the three dimensional L-shaped block with mesh refinement. It can be seen that the proposed integration rule LS1-3D and LS3$n$-3D convergences to the reference solution simultaneously.}
\label{fig:L_shape3d}
\end{figure}

Consider a three dimensional square block with a cubic hole subjected to the surface traction '\textit{t}=1N/mm' as shown in \fref{fig:Lshape}. Only a quarter of the domain is modeled due to symmetry. This results in a three dimensional L-shaped block as shown in \fref{fig:Lshape}. Input parameters used for this analysis are $a$=50mm, $E$=1 MPa and $\nu$=0.3. The reference solution (Strain energy = 382505 MPa) is evaluated using the commercial software Abaqus with a very fine mesh (49211 number of elements) using tetrahedron elements. The result shown in \fref{fig:L_shape3d} concludes that the strain energy converges to the reference solution with the proposed numerical integration rule (i.e. LS1-3D) with mesh refinement. Few of the sample meshes used are shown in \fref{fig:3dmesh}.

\section{Concluding Remarks}
A linearly consistent one point quadrature rule has been proposed to integrate over star convex arbitrary polytopes. The results from the proposed scheme are compared with the linear smoothing scheme. The linear smoothing scheme (denoted as LS3$n$-2D/LS3$n$-3D in the paper) consists in subdividing the poly-element into simplices (triangles or tetrahedra). The linear smoothing scheme is then performed over each triangle. This process requires 3$n$ and 4$n$ integration points per element, where $n$ is the number of sides/face of the element. This significantly reduces the computational effort whilst preserving accuracy and stability. The proposed integration rule yields also preserves optimal convergence rates in both the $L^2$ norm and in the $H^1$ seminorm.

\newpage
\begin{acknowledgements}
St{\'e}phane Bordas thanks the financial support of the European Research Council Starting Independent Research Grant (ERC Stg grant agreement No. 279578) entitled ``Towards real time multiscale simulation of cutting in non-linear materials with applications to surgical simulation and computer guided surgery'' and is also grateful for the support of the Fonds National de la Recherche Luxembourg FNRS-FNR grant INTER/FNRS/15/11019432/EnLightenIt/Bordas. 
\end{acknowledgements}


\bibliographystyle{spphys} 
\bibliography{./bibfile/onepointscheme}


\end{document}